\documentclass[12pt]{article}
\usepackage{amsmath,amsthm,amssymb}

\newtheorem{thm}{Theorem}
\newtheorem{prop}{Proposition}
\newtheorem{lemma}{Lemma}
\newtheorem{defn}{Definition}
\newtheorem{cor}{Corollary}

\newcommand{\ga}{\alpha}
\newcommand{\gb}{\beta} 
\newcommand{\gc}{\gamma}
\newcommand{\gd}{\delta}
\newcommand{\gep}{\epsilon}
\newcommand{\gl}{\lambda}
\newcommand{\go}{\omega}
\newcommand{\gs}{\sigma}

\newcommand{\gO}{\Omega}
\newcommand{\gC}{\Gamma}

\newcommand{\hess}{\nabla^2}
\newcommand{\gra}{\nabla}
\newcommand{\de}{\partial}
\newcommand{\li}{\langle}
\newcommand{\ri}{\rangle}
\newcommand{\rtarw}{\rightarrow}

\newcommand{\bpf}{\begin{proof}}
\newcommand{\epf}{\end{proof}}
\newcommand{\beq}{\begin{equation}}
\newcommand{\eeq}{\end{equation}}
\newcommand{\beqn}{\begin{eqnarray*}}
\newcommand{\eeqn}{\end{eqnarray*}}

\begin{document}
\title{Optimal Curvature Decays on Asymptotically Locally Euclidean Manifolds} 
 \author{Szu-yu Sophie Chen
  \footnote{The author was supported by the Miller Institute for Basic Research in Science.}}
 \date{}
\maketitle

\begin{abstract}
 We present a method in nonlinear elliptic systems to study curvature decays on asymptotically locally Euclidean (ALE) manifolds.
 In particular, we show that scalar flat Kahler ALE metrics of real dimension $n$ are of order $(n-2).$
 The analysis can  also apply to study removable  point singularity problem.  
 \end{abstract}

 Let $(M, g)$ be a Riemannian manifold of dimension $n \geq 4$. We denote the Riemannian curvature, Ricci curvature, the scalar curvature
 and  Weyl curvature by $Rm, Rc, R$ and $W,$ respectively.

 \begin{defn} A complete Riemannian manifold
 $(M^n, g)$ is called  asymptotically locally Euclidean (ALE) of order $\tau > 0$ if there exists
 a compact subset $K \subset M$ such that (each component of) 
  $ M\setminus K$ is diffeomorphic to $(R^n\setminus B_{r_0})/\Gamma,$ 
 where  $\Gamma \subset SO(n)$ is a finite group acting freely on $(R^n\setminus B_{r_0}).$
 Moreover,  under this identification
 $$g_{ij} -\gd_{ij} = O(|x|^{-\tau}),$$
 $$\de^k g_{ij} = O(|x|^{-\tau-k})$$
 as $|x| \rightarrow \infty.$
  $(M^n, g)$ is called ALE of order $0$ if under the above identification,
  $g_{ij} -\gd_{ij} = o(1)$ and  $\de^k g_{ij} = o(|x|^{-k})$ as $|x| \rightarrow \infty.$
\end{defn}

 A classical work by Bando-Kasue-Nakajima \cite{BKN89} asserts that for
 Kahler Ricci-flat ALE metrics there exist coordinates  of order $n,$ and 
 for Ricci-flat metrics,  there exist coordinates of order $n-1.$    
  Cheeger-Tian \cite{CT94} improved the result to  order $n$ and generalized to other cone-like Ricci-flat metrics.   
 For non Ricci-flat metrics, 
 Tian-Viaclovsky \cite{TV05}  proved  in dimension 4 that scalar flat (anti-)self-dual and harmonic metrics  are of order $2- \epsilon.$ Streets \cite{Str07} proved that  scalar flat Bach-flat metrics are of order $2- \epsilon.$

 We study under what conditions, $(M, g)$ is ALE of order $n-2.$ 
 This (optimal) order is the borderline case when the mass is
  finite but not necessary zero; see Bartnik \cite{Bartnik86}. The problem of getting the decay order $n-2$
  is especially delicate, as we need to exclude the possibility of the decay  $|x|^{-(n-2)} \ln |x|;$ this $\ln |x|$ term does 
  not occur in other orders.
  
  To be more precise, we denote by $r$ the distance function
 from a fixed point $o.$ In Theorem~\ref{t:delRm} below, by using the curvature equation
 $ \Delta Rm = \gra \gd Rm + Rm \ast Rm,$
  we  show that the size of $(\gd Rm)_{jkl}= \gra^i R_{ijkl}$ basically controls the decay of the metric.
  However, when $|\gd Rm|$ is asymptotic to $r^{-(n+1)},$  $|Rm|$ is asymptotic to
   $r^{-n} \ln r.$  To eliminate the $\ln r$ term, we show that under a stronger
    condition $|\gra Rc| = O(r^{-(n+1)})$ (this is stronger because $(\gd Rm)_{jkl} =  \gra_k R_{jl}- \gra_l R_{jk}$),
    the metric has the decay rate $n-2.$

 We denote by  $C_s$ the Sobolev constant,  the best constant such that 
 $\|f\|_{L^{\frac{2n}{n-2}}(M)} \leq C_s \|\gra f\|_{L^2(M)} $ for $f \in C^{0,1}_c (M).$ 
 
  \begin{thm}\label{t:delRm} Let $(M, g)$ be ALE of order $0$ with  $C_s < \infty.$  Suppose $|Rm| \in L^{\frac{n}{2}}(M).$

(a) If $\gra^k \gd Rm = O(r^{-(k+s+3)})$  with  $s > \min \{0, \frac{n-6}{2}\},$ then  $(M, g)$ is
   ALE of order $\ga = \min \{n-2, s\}$ when $s \neq n-2,$
        and  $(M, g)$ is ALE of order $\ga < n-2$ when  $s = n-2.$ 
 
 (b) If $\gra^{k}  Rc= O(r^{-(k+n)}),$ then  $(M, g)$ is ALE of order $\ga = n-2.$         
\end{thm}
  In the definition of ALE metrics, if we  only require the decay in lower derivatives, then 
in Theorem~\ref{t:delRm} only smaller $k's$ are needed.

The next result concerns special metrics whose curvatures satisfy a coupled system,
 $$ \left\{ \begin{array} {l}  
 \Delta Rc = Rc \ast Rm\\
  \Delta Rm = \gra \gd Rm + Rm \ast Rm. 
 \end{array} \right. 
 $$
  This allows us to estimate $Rc$ first and then use Theorem~\ref{t:delRm} to obtain the decay rate.

 Let $K$ be a compact subset in $M.$  
 A metric is \emph{harmonic} if $\gd Rm = 0.$ In dimension 4, a metric is \emph{(anti-)self-dual} if $W^- = 0$ ($W^+=0$).
 \begin{thm} \label{t:Kahler} Let $(M, g)$ be a complete noncompact Riemannian manifold with  $R= 0$ in $M \setminus K.$
 Suppose $|Rm| \in L^{\frac{n}{2}}(M)$  and   $C_s< \infty.$ 
 
(a) If $(M, g)$ is Kahler, then $(M, g)$ is ALE of order $n-2$ with finitely many ends.

(b) If  $(M, g)$ is harmonic,  then $(M, g)$ is ALE of order $n-2$ with finitely many ends.

(c) If  $n=4$ and $(M, g)$ is (anti-)self-dual,  then $(M, g)$ is ALE of order $2$ with finitely many ends.
\end{thm}
 In dimension 4, case (a) belongs to case (c); see \cite{LeB88}.

 
 Remark: After the current work has been completed,  the author  was notified  that
   by a different geometry argument, case (c) in Theorem~\ref{t:Kahler} was proved 
 in a recent work
  by Chen-Lebrun-Weber  \cite{CLW08}.   Also, Ricci decay problem for extreme Kahler metrics was considered
   in \cite{CW07}.

 We give the outline of proofs. To fix the notation, 
  we denote by $dV$ the volume element in $M$, and by $d\gs$ the area element of a hypersurface in $M.$ 
  Let $B_r (x)$ and $S_r(x)$ be the  geodesic ball of radius $r$ and the geodesic sphere of radius $r$ centered at $x,$ respectively. 
 When $x$ is at $o$, we simply denote by $B_r$ and $S_r.$ 


  The idea of the proof of Theorem~\ref{t:delRm} is to compare the  size of
   $\int_{M \setminus B_r} |\gra Rm|^2 dV$ (as a function of $r$) and   its derivative $-\int_{S_r} |\gra Rm|^2 d\gs.$ 
  Then by ordinary differential inequality lemma (see Lemma~\ref{l:ode}), we obtain the the decay of $|\gra Rm|$ and hence the decay
   of $|Rm|;$    this is  where $\ln r$ term might occur. Then by a classical result in \cite{BKN89},  there exist coordinates of   the desired order.    To relate above two integrands, we apply   Pohozaev's trick 
     to the system $\Delta Rm = Rm \ast Rm +  \hess Rc$ in the exterior domain to get
    $ \int_{M \setminus B_r} |\gra Rm|^2 dV \leq \frac{2}{n-2} r \int_{S_r} (|\gra_{\nu} Rm|^2- \frac{1}{2} |\gra Rm|^2) d\gs + \text{lower order terms},$   where $\nu$ is the unit outer normal on  $S_r.$ Finally, we apply algebraic inequalities    
   $ \sup_{|v|= 1} |\gra_v Rm|^2 \leq \frac{n}{n+2} |\gra Rm|^2 + C_n |\gd Rm||\gra Rm| + C_n |\gd Rm|^2$ and
    $ \sup_{|v| = 1} |\gra_v Rm|^2 \leq \frac{n-1}{n+1} |\gra Rm|^2$ \\$ + C_n |\gra Rc||\gra Rm| + C_n |\gra Rc|^2$ to obtain the
   sharp comparison between    $|\gra_{\nu} Rm|^2$ and $|\gra Rm|^2$ for cases (a) and (b), respectively.
    It turns out that  the method of using Pohozaev's trick  is flexible and can be applied to 
     more general non-variational elliptic systems; see \cite{Chen08a}.

   For  Theorem~\ref{t:Kahler}, by a  regularity result by Chen-Weber \cite{CW07} 
  and a  work by Tian-Viaclovsky \cite{TV06}, the manifold is ALE of order zero. 
  To improve the order,  we  apply Pohozaev's trick to the equation
    $\Delta Rc = Rm \ast Rc$ to  get $ \int_{M \setminus B_r} |\gra Rc|^2 dV \leq \frac{2}{n-2} r \int_{S_r} (|\gra_{\nu} Rc|^2- \frac{1}{2} |\gra Rc|^2) d\gs + \text{lower order terms}.$
  Then we apply an  inequality 
   $\sup_{|v| = 1}  |\gra_v Rc|^2 \leq$\\ $ \frac{n}{n+2} |\gra Rc|^2,$ which holds for scalar flat
    Kahler and (anti-) self-dual metrics,  to obtain the
   sharp comparison between    $|\gra_{\nu} Rc|^2$ and $|\gra Rc|^2$ and    
  the right fall-off rate of $|Rc|.$ Finally,  applying Theorem~\ref{t:delRm} (b), we find coordinates 
  of the desired order.
 
  The inequality $\sup_{|v| = 1}  |\gra_v Rc|^2 \leq \frac{n}{n+2} |\gra Rc|^2$ is stronger than
  $|\gra |Rc||^2 \leq  \frac{n}{n+2} |\gra Rc|^2.$ The latter is called the Kato inequality.
  A general theory on the Kato inequality  can be found  in \cite{Bran00} and \cite{CGH00} by using the representation theory. 
   We refer the reader to \cite{SSY75}, \cite{BKN89}, \cite{Rade93}, \cite{GL99} for related results  in the literature.

 Now we discuss explicit examples of metrics with the optimal order. 
 The Schwarzchild metric  $(1+ \frac{\mu}{r^{n-2}})^{\frac{4}{n-2}} dx^2$  defined  on $\mathbb{R}^n \setminus \{ 0\}$
 is an example of scalar flat harmonic ALE manifolds of order $n-2.$
 For scalar flat Kahler ALE manifolds, for $n=4,$   Lebrun \cite{LeB88} constructed examples of order $2.$
  In higher dimensions, the problem was studied by Rollin-Singer \cite{RS}  by using the momentum construction. 
Let  $n=2m.$
 \begin{prop}  \cite{RS}  \label{p:kahler} 
There exist scalar-flat Kahler ALE manifolds on  a complex line bundle over $\mathbb{C}P^{m-1}$  of order $n-2$ with $\Gamma= \mathbb{Z}_p$ for all 
 $p \geq 1.$ 
\end{prop}

  We reconstruct the metrics for $n= 2m$ by generalizing the method by LeBrun \cite{LeB88} to higher dimensions.
  The advantage of this approach is that we are able to express the metrics  explicitly. The  metrics  are of the form
  $$ ds^2 = (1+ \ga \rho^{-2(m-1)} + \gb \rho^{-2m})^{-1} d\rho^2 + \rho^2 (1+ \ga \rho^{-2(m-1)} + \gb \rho^{-2m}) h_0 + \rho^2 g_0.$$
   Metrics of this form  were   studied by Simanca \cite{simanca} and Pedersen-Poon \cite{pp}. 
   Explicit examples of scalar flat ALE manifolds were used in the problem of desingularization; see e.g.,
   Arezzo-Pacard \cite{ap06}.

  A by-product of the analysis we used in proving Theorem~\ref{t:delRm} and \ref{t:Kahler} is the following
  removable point singularity result. We show that if the curvatures satisfy a coupled system in a punctured ball, then
  $g$ is  smooth across the origin. 
  The result is a direct analytic consequence,  and  holds on more general orbifolds.  
  \begin{defn} A Riemannian orbifold $(M^n, g)$ is a smooth Riemannian manifold away from finite singular
   points ${x_i}.$ At each point $x_i$, there is a neighborhood $B_r(x_i)$ 
   such that  $B_r(x_i) \setminus \{x_i\}$ is diffeomorphic to a cone on $\mathbb{S}^{n-1}/ \Gamma,$
   where  $\Gamma \subset SO(n)$ is a finite group acting freely on $\mathbb{S}^{n-1}.$ Moreover, if
   $g$ is lifted to $B_1 \setminus \{0\}$ via $\Gamma,$ then $g$ (under a $\Gamma$- equivariant  diffeomorphism) extends 
   to a smooth Riemannian metric on $B_1.$
     
   A $C^0$ orbifold $(M^n, g)$ is a smooth Riemannian manifold away from finite singular
   points with above structure. At each singular point, the lift of $g$ extends to a $C^0$ metric on $B_1.$
  \end{defn}
  \begin{cor} \label{c:remov}
  Let $n \geq 5.$ Let  $(B_1, g)$ be a $C^0$ orbifold. Suppose $(B_1 \setminus \{0\}, g)$ is a smooth n-dimensional 
  Riemannian manifold  satisfying
   $$ \left\{ \begin{array} {l}  
 \Delta Rc = Rc \ast Rm\\
  \Delta Rm = \gra \gd Rm + Rm \ast Rm 
 \end{array} \right. 
 $$ on $B_1\setminus \{0\}$ with $|Rm| \in L^{\frac{n}{2}} (B_1)$ and $C_s < \infty.$ Then
  $g$  extends to a smooth orbifold metric on $B_1.$ 
   \end{cor}
 Remark:    Corollary~\ref{c:remov}
  can be viewed as a Riemannian analogue of Sibner's \cite{Sibner85}  result for 
  coupled Yang-Mills systems. The argument we used here was applied to  extreme Kahler metrics  by  Chen-Weber \cite{CW07}.
 
  Manifolds satisfying the system in Corollary~\ref{c:remov} include constant 
 scalar curvature Kahler metrics and harmonic metrics. 
   
 At the end of the introduction, we state the geometrical Pohozaev's identity for an  independent interest. Let $\Gamma$ be the Christoffel symbol.
 \begin{prop}\label{p:poh} Let $\Omega$ be a domain in $M$ and $T$ be a tensor. Suppose the coordinate vector $X= X^i \gra_i$ 
 is defined in $\Omega.$  Then
  \beqn
  \int_{\gO} \li \Delta T, X^i \gra_i T \ri dV&=& \int_{\gO} \left( \frac{n-2}{2} |\gra T|^2 + \gra T \ast \gra T \ast \gC \ast X +   \gra T \ast  T \ast Rm \ast X \right) dV\\
  &+& \int_{\de \gO}  \li \gra_{\nu} T, X^i \gra_i T\ri d\gs-\int_{\de \gO} \frac{1}{2} \li X, \nu \ri |\gra T|^2 d\gs,
 \eeqn where $\nu$ is the unit outward normal on $\de \Omega.$ 
 \end{prop}
 The Pohozaev's identity was used in the literature for  the Yamabe equation \cite{KMS07}. 

 The organization of the paper is as follows.  In Section~\ref{s:prelim}, we collect  some standard regularity results and 
 review background materials.
 In Section~\ref{s:kahler},  we construct scalar-flat Kahler metrics. We prove algbraic inequalities in Section~\ref{s:Kato}, and 
 Theorem~\ref{t:delRm}   and \ref{t:Kahler} in Sections~\ref{s:thm}. Finally, we prove Corollary~\ref{c:remov} in Section~\ref{s:cor}

   \textbf{Acknowledgments:}  The author would like to thank Alice Chang and Paul Yang for stimulating  discussions which initiated the present work.      
     She appreciates  Jeff Viaclovsky for helpful discussions and especially for explaining  his work to her.  
      The author is  grateful to Rick Schoen for helpful suggestions, which result in Theorem~\ref{t:delRm}.
      The author thanks Claude LeBrun, Frank Pacard and Jian Song for useful discussions, and Gang Tian for interests.
      Finally, the author  thanks Claude LeBrun for pointing out the article  \cite{CLW08} to her and Peter Petersen
      for providing a reference.
      
\section{Preliminaries} \label{s:prelim}
  We collect some standard results in  elliptic regularity theory and ordinary differential equations. 
 Then we review Kahler, harmonic and (anti-)self-dual metrics.
   
   Let $D_r$ be  the complement of the geodesic ball of radius $r.$
 \begin{lemma} \cite{BKN89} \label{l:homo}  Let $p > 1.$ Suppose $u, f \geq 0$ and  $u \in L^p (D_r)$ satisfies  $\Delta u \geq -  f u$ 
   in $D_r.$
  
  (a)  Suppose $f \in L^q (D_r)$ satisfies $\int_{D_r} f^q dV \leq C_0 r^{-(2q-n)}$ for some $q > \frac{n}{2}.$  Then 
         $\sup_{D_{2r}} u \leq C r^{-\frac{n}{p}} \|u\|_{L^p(D_r)},$ where $C = C(n, p, C_0, C_s)$
   
  (b)  Assume the conditions in (a) and in addition that $f \in L^{\frac{n}{2}} (D_r)$ and $Vol B_r(x) \leq C r^n.$ Then
           $u = O(|x|^{- \ga}) \; \text{for all}\;\, \ga < n-2 \;\text{as} \;\, |x| \rtarw \infty.$
 \end{lemma}

Next lemma is an analogue of Lemma~\ref{l:homo} for  inhomogeneous equations.
The proof follows by an argument similar to  \cite{BKN89}, section 4.
 Let $\gc = \frac{n}{n-2}.$
 \begin{lemma} \label{l:inhomo}
   Let  $u, f, h \geq 0.$    Suppose $u$ satisfies $\Delta u \geq -  f u- h$ in $D_r$
   and  $Vol B_r (x) \leq C_1 r^n.$     
     
  (a)  Suppose $f \in L^q (D_r)$ satisfies $\int_{D_r} f^q dV \leq C_0 r^{-(2q-n)}$ for some $q > \frac{n}{2}.$ 
       Assume  that  $u \in L^p(D_r)$ for some $p> 1$ and $h \in L^q (D_r).$
        Then 
         $\sup_{D_{2r}} u \leq C r^{-\frac{n}{p}} \|u\|_{L^p(D_r)}+C r^{\frac{2q-n}{q}} \|h\|_{L^q(D_r)},$
        where $C=C(n, p, C_0,  C_1, C_s).$

  (b)   Let $q_1 = \frac{p \gc}{p \gc - p + 1}= \frac{pn}{2p+n-2}$ and $p  > 1.$  
       Suppose  $h \in L^{q_1}(D_r)$ satisfies $\|h\|_{L^{q_1}(D_r)} < C r^{-\gd_1}$ for some $\gd_1 > 0.$ 
        Assume $u \in L^{p \gc} (D_r)$ and  $f \in L^{\frac{n}{2}} (D_r).$
       Then there exists $\gep_1 > 0$ such that
        if  $\|f\|_{L^\frac{n}{2}(D_r)} < \epsilon_1,$ then
       $\int_{D_r} u^{p\gc} dV= O(r^{-\epsilon}) \; \text{for some} \; \epsilon > 0 \;\text{as} \;\, r \rtarw \infty.$

       Moreover, under the same assumptions as above but $u \in L^p (D_r)$ for $p > 1,$ it holds
       $$\left( \int_{D_{2r}}  u^{p \gc} dV \right)^{\frac{1}{\gc}} \leq C r^{-2} \int_{D_r} u^p dV+ C \left(\int_{D_r} h^{q_1} dV\right)^{\frac{p}{q_1}} .$$  

  (c)  Assume the conditions in (a) and in addition that $f \in L^{\frac{n}{2}} (D_r)$ and $\|h\|_{L^q(D_r)} \leq Cr^{-(s+4- \frac{n}{q})}.$
        Then
           $u = O(|x|^{- \ga}) \; \text{for all}\;\, \ga < \min \{n-2, s+2\} \;\text{as} \;\, |x| \rtarw \infty.$
 \end{lemma}
 \bpf (a) is by standard elliptic regularity; see \cite{Morrey}
  
  (b) 
   Let $0 \leq \phi \leq 1$ be a cutoff function such that  $\phi = 0$ on $B_r \cup  D_{2r'}$
    and $\phi = 1$  on $B_{r'}\setminus B_{2r}$ with $|\gra \phi| \leq C r'^{-1}$ on $B_{2 r'} \setminus B_{r'}$
    and  $|\gra \phi| \leq C r^{-1}$ on $B_{2r} \setminus B_{r}.$ $r'$ will be chosen large.
   Applying $\phi^2 u^{p-1}$ to the equation and following the computations
   in \cite{BKN89} p332, we get
        $\left( \int |\phi u^{\frac{p}{2}}|^{2 \gc} dV \right)^{1/\gc} 
        \leq C \int \left(f \phi^2 u^p + \phi^2 u^{p-1} h + |\gra \phi|^2 u^p\right) dV.$
   Since $\|f\|_{L^\frac{n}{2}}$ is small, we can absorb  the term $ \int f \phi^2 u^p dV$ to the left and get
        \beqn\left( \int |\phi u^{\frac{p}{2}}|^{2 \gc} dV \right)^{1/\gc} 
                &\leq& \left(\int \phi^{2 q'} u^{p \gc} dV\right)^{1/q'}
        \left( \int_{supp \phi} h^{q_1} dV\right)^{1/q_1} +  C\int |\gra \phi|^2 u^p dV,
        \eeqn where $q' = \frac{p \gc}{p-1}$ satisfies $\frac{1}{q_1} + \frac{1}{q'} = 1.$
        Note that $\frac{q'}{\gc} = \frac{p}{p-1}.$ We use Yang's inequality $ab \leq \frac{(a/\gep)^p}{p} + \frac{p-1}{p}(\gep b)^{\frac{p}{p-1}}$
        to lift the power $\frac{1}{q'}$ to $\frac{1}{\gc}$
        and absorb  the term to the left to obtain
         \beq \label{i:inhomo-reg}
         \left( \int |\phi u^{\frac{p}{2}}|^{2 \gc} dV \right)^{1/\gc} 
        \leq C \left( \int_{supp \phi} h^{q_1} dV\right)^{p/q_1} +  C\int |\gra \phi|^2 u^p dV.
        \eeq
        Using Holder inequality for the second term, we get 
         \beq \label{i:iterate}
         \left( \int |\phi u^{\frac{p}{2}}|^{2 \gc} dV \right)^{\frac{1}{\gc}} 
        \leq C \left( \int_{supp \phi} h^{q_1} dV\right)^{\frac{p}{q_1}} +  
        C\left(\int |\gra \phi|^n dV \right)^{\frac{2}{n}} \left(\int_{supp |\gra \phi|} u^{p\gc} dV\right)^{\frac{1}{\gc}}.
        \eeq
        Note that $\int_{D_{r'} \setminus D_{2r'}} u^{p\gc} dV \rtarw 0$
         because $u \in L^{p \gc}.$
       Letting $r' \rtarw \infty$, we have
        $\int_{D_{2r}}  u^{p \gc} dV \leq C r^{- \gd_1 p \gc} + C  \int_{D_r \setminus D_{2r}}  u^{p \gc} dV.$
        Therefore, $\int_{D_r}  u^{p \gc} dV= O(r^{-\gep})$ for some $\gep > 0.$
       
         For the second claim, we go back to (\ref{i:inhomo-reg}).
         Letting $r' \rtarw \infty$ gives the inequality.
        
  (c) Let $\bar \ga= \sup \{\ga : u = O(|x|^{-\ga})\}.$ Suppose $\bar \ga < \min \{n-2, s+2\}.$
        Let $\tilde p \, ( = p \gc ) > \frac{n}{\bar \ga} > \gc$ be close to $\frac{n}{\bar \ga}.$
       Then $u \in L^{\tilde p}.$
      We have $q_1 = \frac{ p \gc}{ p \gc -  p + 1}< \frac{n}{2}< q.$ 
      Therefore, 
      $\|h\|_{L^{q_1}  (D_r)} \leq Cr^{-\frac{n}{q}+\frac{n}{q_1}} \|h\|_{L^{q} (D_r)}  \leq Cr^{-(s+4)+ \frac{n}{q_1}}.$
      Since $\frac{n}{q_1}$ is close to $\frac{n (\gc-1)+ \gc \bar \ga}{\gc}= 2+ \bar \ga < s+4,$
      by (b) we get $\int_{D_r} u^{\tilde p} = O(r^{-\epsilon})$ for some $\gep> 0.$ By (a), it follows 
       $\sup_{2 D_r} u \leq C r^{-n/\tilde p} \|u\|_{L^{\tilde p}(D_r)}+ Cr^{-(s+2)}= O(r^{-\frac{n+ \gep}{\tilde p}}+ r^{-(s+2)}).$
      This gives a contradiction as $\frac{n+ \gep}{\tilde p}, s+2 > \bar \ga.$
 \epf
We  turn to consider the punctured ball $B_1 \setminus \{0\}$ and
  study the blow-up rate of  solutions near the origin. 
  Let $A_r = \{x: \frac{r}{2} < |x| < \frac{3r}{2} \}$ be an annulus.
 \begin{lemma} \label{l:sing-homo} Let $n \geq 5.$
 Suppose that $u, f \geq 0$ are in $L^{\frac{n}{2}} (B_1)$ and $u$ satisfies  $\Delta u \geq -  f u$  in $B_1 \setminus \{0\}.$ 
 Assume $Vol B_r(x) \leq C r^n$ and in addition $\int_{A_r} f^q dV \leq C r^{-(2q-n)}$ for some $q > \frac{n}{2}.$ 
  Then there exists $\epsilon_0 > 0$ such that if 
  $(\int_{B_{r_0}} f^{\frac{n}{2}} dV)^{\frac{2}{n}} \leq \epsilon_0,$ then  
    $ \int_{B_r} u^{\frac{n}{2}} dV = O(r^{\epsilon}) \; \text{for some} \; \epsilon> 0$ as $r \rtarw 0.$
  
   Moreover, 
   $u = O(|x|^{- \ga}) \; \text{for all}\; \ga > 0 \;\text{as} \;\, |x| \rtarw 0.$
 \end{lemma}
 \bpf We  again follow the computations   in \cite{BKN89} p332.
   Let $\bar q= p \gc > \gc.$  Suppose $u \in L^{\bar q} (B_1).$
   Let $0 \leq \phi \leq 1$ be a cutoff function such that  $\phi = 0$ on $B_{r'} \cup (B_1 \setminus B_{2r})$
    and $\phi = 1$  on $B_{r}\setminus B_{2r'}$ with $|\gra \phi| \leq C r'^{-1}$ on $B_{2 r'} \setminus B_{r'}$
    and  $|\gra \phi| \leq C r^{-1}$ on $B_{2r} \setminus B_{r}.$ $r'$ is chosen to be small. Then by   
     $(\int_{B_{r_0}} f^{n/2} dV)^{2/n} \leq \epsilon_0,$ we have (by (\ref{i:iterate}))
    $$\left( \int |\phi u^{\frac{p}{2}}|^{2 \gc} dV \right)^{\frac{1}{\gc}} 
        \leq  
        C\left(\int |\gra \phi|^n dV \right)^{\frac{2}{n}} \left(\int_{supp |\gra \phi|} u^{p\gc} dV\right)^{\frac{1}{\gc}}
        \leq C \left(\int_{supp |\gra \phi|} u^{p\gc} dV\right)^{\frac{1}{\gc}}.$$
  Letting $r' \rtarw 0$ (and noting $u \in L^{\bar q}$), we get
   $\int_{B_r} u^{\bar q} dV \leq C \int_{B_{2r}\setminus B_r} u^{\bar q} dV.$
   Therefore, $\int_{B_r} u^{\bar q} dV = O(r^{\epsilon})$ for some $\epsilon > 0.$
    Now let $\bar q= n/2.$ We obtain the first part of the lemma.
    
     By standard regularity, 
     $\sup_{|x|= r} u(x) \leq C r^{- n/p} (\int_{A_r} u^p dV)^{1/p}$ if $u \in L^p(B_1).$ Since  $u \in L^{\frac{n}{2}} (B_1),$
     we have $u = O(|x|^{-2})$ near the origin. Let $\bar \ga = \inf \{\ga > 0: u=O(|x|^{-\ga})\}.$
     If $\bar \ga > 0,$  then $\gc = \frac{n}{n-2} < \frac{n}{2} \leq \frac{n}{\bar \ga}.$
     Let $\gc < \bar q < \frac{n}{\bar \ga}$ be close to $\frac{n}{\bar \ga}.$ Then $u \in L^{\bar q} (B_1).$ By the computation
      in the previous paragraph, we get $\int_{B_r} u^{\bar q} dV = O(r^{\epsilon})$ for  some $\epsilon > 0.$
      Hence, $\sup_{|x|=r} u(x) \leq C r^{- n/\bar q} (\int_{A_r} u^{\bar q} dV)^{1/{\bar q}} = O(r^{- \frac{n}{\bar q} + \frac{\epsilon}{\bar q}}).$
     Since $\frac{n}{\bar q} - \frac{\epsilon}{\bar q} < \bar \ga,$ we get a contradiction.
\epf

 \begin{lemma} \label{l:sing-inhomo} Let $n \geq 5.$
  Suppose $h \geq 0$ satisfies $\|h\|_{L^{q}(A_r)} \leq C r^{-2-\gd+ \frac{n}{q}}$ for some $q > \frac{n}{2}$
   and for all $\gd> 0.$
  Suppose $u, f \geq 0$ are in $L^{\frac{n}{2}} (B_1)$ and $u$ satisfies  $\Delta u \geq -  f u- h$  in $B_1 \setminus \{0\}.$ 
   Assume $Vol B_r (x) \leq C r^n$ and in addition $\int_{A_r} f^q dV \leq C r^{-(2q-n)}.$      Then there exists $\epsilon_0 > 0$ such that if 
  $(\int_{B_{r_0}} f^{\frac{n}{2}} dV)^{\frac{2}{n}} \leq \epsilon_0,$ then 
   $u = O(|x|^{- \ga})\;  \text{for all} \;\, \ga > 0 \;\text{as} \;\, |x| \rtarw 0.$
 \end{lemma}
 \bpf The proof is similar to Lemma~\ref{l:inhomo}.
  Let $\bar q= p \gc > \gc.$  Suppose $u \in L^{\bar q} (B_1).$
   Let $0 \leq \phi \leq 1$ be a cutoff function as in Lemma~\ref{l:sing-homo}.   By   
     $(\int_{B_{r_0}} f^{n/2} dV)^{2/n} \leq \epsilon_0,$ the formula (\ref{i:iterate}) holds.
  Letting $r' \rtarw 0$ in (\ref{i:iterate}) and noting that $u \in L^{\bar q}$, we get
   \beq \label{i:inhomo}
   \left(\int_{B_r} u^{\bar q} dV \right)^{1/\gc} \leq C \left( \int_{B_{2r}} h^{q_1} dV\right)^{p/q_1}
   + C \left(\int_{B_{2r}\setminus B_r} u^{\bar q} dV\right)^{1/\gc}.
   \eeq

  By standard regularity, 
     $\sup_{|x| = r} u(x) \leq C r^{- n/p} (\int_{A_r} u^p dV)^{1/p}+ C r^{\frac{2q-n}{q}} \|h\|_{L^q(A_r)}$
      if $u \in L^p(B_1).$ 
     Since  $u \in L^{n/2} (B_1),$ we get    $u = O(|x|^{-2})$ near the origin. Let $\bar \ga = \inf \{\ga > 0: u=O(|x|^{-\ga})\}.$
     If $\bar \ga > 0,$  then $\gc < \frac{n}{2} \leq \frac{n}{\bar \ga}.$
     Let $\gc <\bar q  < \frac{n}{\bar \ga}$ be close to $\frac{n}{\bar \ga}.$ Then $u \in L^{\bar q} (B_1).$ 
     Since $q_1 = \frac{p \gc}{p \gc - p + 1}< n/2 < q,$ we have
      $\|h\|_{L^{q_1}(A_r)} \leq C r^{- \frac{n}{q} + \frac{n}{q_1}} \|h\|_{L^q (A_r)} \leq C r^{-2 - \gd + \frac{n}{q_1}}.$
      Since $\bar q$ is close to $\frac{n}{\bar \ga},$ then $\frac{n}{q_1}$ is close to $2 + \bar \ga.$  We obtain
       $-2 - \gd + \frac{n}{q_1} > 0.$ 
      As a result, $\|h\|_{L^{q_1}(A_r)} \leq C r^{\epsilon_1}$ for some $\epsilon_1 > 0.$
      By (\ref{i:inhomo}), 
      $\int_{B_r} u^{\bar q} dV \leq C r^{\epsilon_1 q} + C \int_{B_{2r}\setminus B_r} u^{\bar q} dV.$
      Therefore, $\int_{B_r} u^{\bar q} dV= O(r^{\epsilon_2})$ for some $\epsilon_2 > 0.$ 
        Hence, $\sup_{|x| = r} u(x) \leq C r^{- n/\bar q} (\int_{A_r} u^{\bar q} dV)^{1/\bar q} 
        + C r^{-\gd}= O(r^{- \frac{n}{\bar q} + \frac{\epsilon_2}{\bar q}}+ r^{-\gd}).$
     Since $\frac{n}{\bar q} - \frac{\epsilon_2}{\bar q}, \gd  < \bar \ga,$ we get a contradiction.
 \epf

 Now we  apply  Lemma~\ref{l:inhomo} to the curvature equation derived by Lichnerowicz in 1950:
  \beq \label{e:curv}
\Delta Rm = \gra \gd Rm + Rm \ast Rm. 
\eeq

 \begin{prop} \label{t:reg-curv} Let $n \geq 4.$ Suppose $\gra^k \gd Rm = O(r^{-(k+s+3)})$ and $Vol B_r (x) \leq C_1 r^n.$
  Then there exists $\epsilon > 0$ such that
  if $\|Rm\|_{L^{n/2}(D_r)} < \epsilon,$ 
  $$\sup_{D_{2r}} |\gra^k Rm| \leq C r^{-(k+ 2)} \|Rm\|_{L^{\frac{n}{2}}(D_r)} + C r^{-(k+s+2)},$$
  where $C=C(n, k, C_1, C_s).$
 \end{prop}
 \bpf  
  By (\ref{e:curv}), $\Delta |Rm| \geq -C  |Rm| |Rm|-C |\gra \gd Rm|.$
  By Lemma~\ref{l:inhomo}, the case $k= 0$ holds.
  Inductively we have
  \beq \label{e:curv-deriv}\Delta \gra^i Rm = Rm \ast \gra^i Rm + ( \gra^{i+1} \gd Rm+ \sum_{l=1}^{i-1} \gra^l Rm \ast \gra^{i-l} Rm), \eeq
  \beq\label{i:curv-deriv}\Delta |\gra^i Rm| \geq -C |Rm| |\gra^i Rm| - C ( |\gra^{i+1} \gd Rm|+ \sum_{l=1}^{i-1} |\gra^l Rm| |\gra^{i-l} Rm|).\eeq
  Assume for $0 \leq i \leq k-1$ we have pointwise bounds. Let $\phi$ be a cutoff function such that $\phi= 1$ on $B_{2r}\setminus B_r$ and
   $\phi= 0$ on $B_{r/2} \cup D_{5r/2}$ with $|\gra \phi| \leq C r^{-1}.$ Multiplying (\ref{e:curv-deriv})  with $i = k-1$  by $\phi^2 \gra^{k-1} Rm,$
    we get (after applying Schwarz inequality for the cross term)
    \beqn
    \int |\phi \gra^k Rm|^2 \leq \int \phi^2 |Rm| |\gra^{k-1}Rm|^2 + \int |\gra \phi|^2 |\gra^{k-1}Rm|^2 \\
    + \int\phi^2 |\gra^{k-1} Rm| (|\gra^k \gd Rm|+ \sum_{l=1}^{k-2} |\gra^l Rm| |\gra^{k-1-l} Rm|). 
     \eeqn By  induction hypothesis and the assumption on $\gra^{k+1} \gd Rm,$ 
     we have 
     $$\int_{B_{2r}\setminus B_r} |\gra^k Rm|^2\leq Cr^{-2(k+2)+ n} (\int_{D_r} |Rm|^{\frac{n}{2}})^{\frac{4}{n}}
     + C r^{-2(k+s+2)+n}.$$
      Hence,  applying  Lemma~\ref{l:inhomo} (a) 
    to the equation (\ref{i:curv-deriv})  with  $i=k$ and $p=2$  (which is also true if we consider the domain $B_{2r}\setminus B_r$ instead of $D_r$), we get 
      $$\sup_{|x|= \frac{3r}{2}}|\gra^k Rm| \leq C r^{-\frac{n}{2}} \|\gra^k Rm\|_{L^2(B_{2r}\setminus B_r)}+ C r^{-(k+s+2)}
       \leq Cr^{-(k+2)} \|Rm\|_{L^{\frac{n}{2}}(D_r)} + C r^{-(k+s+2)}.$$
      Hence, $\sup_{D_{2r}} |\gra^k Rm| \leq C r^{-(k+ 2)} \|Rm\|_{L^{\frac{n}{2}}(D_r)} + C r^{-(k+s+2)}.$
 \epf
 
  The following ordinary differential inequality lemma will play an important role in
 proving  Theorem~\ref{t:delRm} and \ref{t:Kahler}.
 \begin{lemma} \label{l:ode}
  Suppose $f(r) \geq 0$ satisfies $f(r) \leq - \frac{r}{a} f'(r) + C_0 r^{-b}$ for some $a, b > 0.$
  
  (a) $a \neq b.$ Then there exists $C_1> 0$ such that
      $f(r) \leq C_1 r^{-a} + \frac{a \,C_0}{a-b} r^{-b}.$ Therefore,  $$f(r) = O(r^{-\min \{a, b\} })
      \quad \text{as}\; r \rtarw \infty.$$
  
  (b) $a = b.$ Then there exists $C_1> 0$ such that
       $f(r) \leq C_1 r^{-a} + a \,C_0 r^{-a}\ln r.$ Therefore,  $$f(r) = O(r^{-a} \ln r)\quad \text{as} \; r \rtarw \infty.$$
 \end{lemma}
 \bpf By the  equation, we have
 $(f(r) r^a)' = r^a f'(r) + a r^{a-1} f(r) \leq a C_0 r^{-b-1+a}.$
 Therefore, $f(r) r^a \leq C_1 + a C_0 \int s^{-b-1+a} ds.$
 If $a \neq b,$  we obtain $f(r) r^a \leq C_1 + \frac{a C_0}{a-b} r^{-b+a}.$
 If $a = b,$ then we get $f(r) r^a \leq C_1 + a C_0 \ln r.$ Multiplying by $r^{-a}$ gives the result.
 \epf

 We review some basic facts about  Kahler, harmonic and (anti-)self-dual metrics.
 
 Let $(M, g)$ be a Kahler manifold with coordinates $Z^1, ..., Z^m$ and an associated complex structure $J.$ 
 Then $\de_{Z_i} = \frac{1}{\sqrt 2} (e_i - \sqrt{-1}\, e_{Ji})$ and $\{e_1,\cdots, e_m, e_{J1}, \cdots, e_{Jm}\}$ is a basis in real coordinates.
  We give the relation of curvatures between real and complex coordinates; see \cite{Joyce}.
  For $1 \leq i,j \leq m,$ we have that $R_{ij} = R_{Ji Jj}$ is symmetric and $R_{i Jj}= - R_{Ji\, j}$ is skew-symmetric
 (but note that $R_{i Jj} = R_{Jj \,i}$). Therefore, $R = 2\, tr R_{ij}.$ The Ricci curvatures in complex coordinates are
  $R_{Z_i \bar Z_j} = R_{ij}+ \sqrt{-1} R_{i Jj},$ which satisfy the Bianchi  relation  $R_{Z_i \bar Z_j, Z_k}=R_{Z_k \bar Z_j, Z_i}.$
   Hence, $$R_{Z_i \bar Z_j, Z_k \bar Z_k}= R_{Z_k \bar Z_j, Z_i \bar Z_k}= R_{Z_k \bar Z_k, \bar Z_j  Z_i}
    + Rc(Z, \bar Z)\ast R(Z, \bar Z, Z, \bar Z).$$
 If $R$ is constant, then $R_{Z_k \bar Z_k, \bar Z_j  Z_i}
 = \frac{1}{2} R_{,\bar Z_j  Z_i}= 0.$  In real coordinates, we get $\Delta Rc = Rc \ast Rm.$
 (To be more precise, $\Delta R_{ij}$ comes from the real part and $\Delta  R_{i Jj}$ comes from the imaginary part.)
 Moreover, 
    \beq \label{e:K-bianchi}
   R_{Z_i \bar Z_j, Z_k} = \frac{1}{\sqrt 2} (R_{ij,k}+ R_{i Jj, Jk}+ \sqrt{-1} R_{i Jj, k}- \sqrt{-1} R_{ij, Jk}).
   \eeq
 
  Let $(M, g)$ be harmonic, i.e. $R_{ijkl, i} = 0.$ By Bianchi identity this gives $R_{jl,k}= R_{jk,l}$ and  
  $R$ is constant. Therefore,   $R_{ij, kk}= R_{ik, jk}= R_{kk, ij} + Rc \ast Rm = Rc \ast Rm.$
   Hence, we have  $\Delta Rc = Rc \ast Rm.$ Examples of harmonic metrics include Einstein, constant scalar curvature locally 
   conformally flat and parallel Ricci curvature (i.e., $\gra Rc= 0$) metrics.
 
  We turn to (anti-)self-dual metrics. In dimension 4,  two forms can be decomposed into $\Lambda^2 = \Lambda^2_+ \oplus \Lambda^2_-,$
  where  $\Lambda^2_+$($\Lambda^2_-$) is the eigenspace of eigenvalue $+1$($-1$) of the Hodge operator. Correspondingly, $W$ is 
  decomposed into $W^+$ and $W^-$ called self-dual and anti-self-dual part.  If $W^-= 0$ ($W^+=0$), we call it a self-dual (anti-self-dual) metric.  
   It can be shown \cite{Besse87} that if $W^- = 0$ (or $W^+= 0$), then $B_{ij} := W_{ikjl, lk} + \frac{1}{2} R_{kl} W_{ikjl}= \Delta R_{ij} + \hess R + Rc \ast Rm = 0.$
  Hence, if $R$ is constant, we have $\Delta Rc = Rc \ast Rm.$
 
 Since every metric satisfies (\ref{e:curv}), the above examples  satisfy the following  system:
 \beq \left\{ \begin{array} {l} \label{e:Rc} 
 \Delta Rc = Rc \ast Rm\\
  \Delta Rm = \gra \gd Rm + Rm \ast Rm= \hess Rc + Rm \ast Rm. 
 \end{array} \right. 
\eeq
 The  $\epsilon$-regularity result for (\ref{e:Rc}) was proved by Chen-Weber \cite{CW07}.  
  \begin{thm} \cite{CW07} Let $n \geq 4.$ Suppose $(M, g)$ satisfies (\ref{e:Rc}). Let $q > n/2$ and $k$ be a nonnegative integer. Then
     there exists $\epsilon > 0$ such that if $\int_{B_r} |Rm|^{n/2} < \epsilon,$ we have
     \begin{align}
     (\int_{B_{r/2}} |\gra^k Rc|^q)^{\frac{1}{q}} \leq C_1 r^{-k-2+\frac{n}{q}} (\int_{B_r} |Rc|^{\frac{n}{2}})^{\frac{2}{n}},\label{Rc-reg}\\
      (\int_{B_{r/2}} |\gra^k Rm|^q)^{\frac{1}{q}} \leq C_1 r^{-k-2+\frac{n}{q}} (\int_{B_r} |Rm|^{\frac{n}{2}})^{\frac{2}{n}}, \label{Rm-reg}
     \end{align} 
     \beq \label{CW}   
     \sup_{B_{r/2}} |\gra^kRm| \leq C_2 r^{-k-2}   \|Rm\|_{L^{\frac{n}{2}}(B_r)} ,
     \eeq where $C_1= C_1(C_s,k, q, n)$ and  $C_2= C_2(C_s,k, n).$
  \end{thm}


\section{Scalar flat Kahler ALE manifolds} \label{s:kahler}

We construct  one parameter family of scalar flat Kahler metrics on the blow up of  $\mathbb{C}^m$
 modulo $\mathbb{Z}_p,$ which is a complex line bundle over $\mathbb{C}P^{m-1}.$ 
 The construction  bases on an idea  by LeBrun \cite{LeB88} for $n=4$. Let $n= 2m.$   
 \bpf[Proof of Proposition~\ref{p:kahler}]
 Let  $\phi (u)$ be  a Kahler potential, where $u = \sum_{i= 1}^m  |z_i|^2.$ Thus,
  $\go = \frac{1}{2} \sqrt {-1} \de \bar{\de} \phi(u)=\frac{1}{2}  \sqrt {-1} g_{i \bar j} dz^i \wedge dz^{\bar j}$ and
  $g_{i \bar j}= \phi'(u) \gd_{ij} + \phi''(u) z_{\bar i} z_j.$
  Let $\psi = \ln \det g_{i \bar j}.$ We have $\psi  = \ln \left( \phi'(u)^{m-1} (\phi'(u) + \phi''(u) u)\right).$ 
 Then the scalar-flat equation becomes
  $$
  0= g^{i \bar j} \psi_{i \bar j} =  g^{i \bar j} (\psi'(u) \gd_{ij} + \psi''(u) z_{\bar i} z_j)
    = \phi'(u)^{m-1} (m\psi'(u) + \frac{\phi'(u)\psi''(u) - \phi''(u) \psi'(u)}{\phi'(u) + \phi''(u) u}  u).
  $$
 In other words, we plan to solve 
 \beq \label{e:ode}
 m\phi'(u)\psi'(u) + (m-1) \phi''(u) \psi'(u) u + \phi'(u) \psi''(u) u= 0
 \eeq
  together with $\psi  = \ln \left( \phi'(u)^{m-1} (\phi'(u) + \phi''(u) u)\right).$

We  observe that
 $$ [\phi'(u)^{m-1} u^m \psi'(u)]' = u^{m-1} \phi'(u)^{m-2} [m\phi'(u)\psi'(u) + (m-1) \phi''(u) \psi'(u) u + \phi'(u) \psi''(u) u].$$ 
 Hence, (\ref{e:ode}) is equivalent to
 $$
  \phi'(u)^{m-1}  \psi'(u) u^m = \ga $$ with
 $ \psi  = \ln \left(\phi'(u)^{m-1} (\phi'(u) + \phi''(u) u)\right),$ 
  where $\ga$ is a constant.
 Let $y(u) = u \phi'(u).$ Then $\psi = \ln \frac{y^{m-1} y'}{u^{m-1}}$
 and $\psi' = \frac{m-1}{y} y' + y'' (y')^{-1} - \frac{m-1}{u}.$ The equation becomes
 $$\ga y' = (m-1) u y^{m-2}  (y')^2 + y^{m-1} u y'' - (m-1) y^{m-1} y' = (u y^{m-1} y' - y^m)'.$$ 
  Thus,
  \beq \label{e:sode}
  y^m + \ga y + \gb = u y^{m-1} y'.
  \eeq
  Formally, we get
 \beq \label{e:u}
  \ln u = \int \frac{y^{m-1} dy}{y^m + \ga y + \gb}.
  \eeq

  To express the metric, note that the standard metric on $\mathbb{S}^{2m-1}$ is decomposed into 
  $g_0$, the standard metric on $\mathbb{C}P^{m-1},$ and $h_0,$ the metric along  the  fiber (the Hopf map). 
  Hence, $ds^2 = y' [(d \sqrt u)^2 + u h] + y g_0.$ Let $y = \rho^2.$ By (\ref{e:sode}), we have
  $y' (d \sqrt u)^2 = y' (4u)^{-1} du^2 = \frac{y}{u y'} d\rho^2=  (1+ \ga \rho^{-2(m-1)} + \gb \rho^{-2m}) d\rho^2$ and $y' u = \rho^2 (1+ \ga \rho^{-2(m-1)} + \gb \rho^{-2m}).$
  Thus, 
  \beq \label{metric}
  ds^2 = (1+ \ga \rho^{-2(m-1)} + \gb \rho^{-2m})^{-1} d\rho^2 + \rho^2 (1+ \ga \rho^{-2(m-1)} + \gb \rho^{-2m}) h_0 + \rho^2 g_0.
  \eeq
 
   We exhibit  some special solutions which give complete metrics.
 Let $\ga = (p-m) a^{2(m-1)}$ and $\gb= (m-1-p) a^{2m}.$  The metric in (\ref{metric}) is defined for $\rho^2 \geq a^2.$  
 When $m=2,$ they coincide with the metrics in  \cite{LeB88}.
  We have 
  $$y^m + \ga y + \gb  = (y- a^2) (y^{m-1}+ a^2 y^{m-2} + \cdots + a^{2(m-2)}y + (p+1 -m) a^{2(m-1)}).$$
   Therefore, $\frac{y^{m-1}}{y^m + \ga y + \gb} > 0$ and  (\ref{e:u}) is invertible for $y \in (a^2, \infty)$ (and  $u \in (0, \infty)$), 
   which gives implicitly  $\phi'(u)= u^{-1} y (u)$ and  $\phi(u).$
 We show that the metric is complete  on the blow up of $\mathbb{C}^m$ module $\mathbb{Z}_p.$ 
   Denote $h_0 = d\theta^2.$ Introducing the coordinates $\hat r^2 = \rho^2- a^2$ and $\hat \theta= p \theta.$
  We get $d \rho^2 = \frac{\hat r^2}{\hat r^2 + a^2} d \hat r^2$ and 
  \beqn
  ds^2|_{\text{fiber}} &=& (1- (\frac{a}{\rho})^2)^{-1} (1+ (\frac{a}{\rho})^2 + \cdots + (\frac{a}{\rho})^{2(m-2)} + (p+1 -m) (\frac{a}{\rho})^{2(m-1)})^{-1} d\rho^2\\
                       & & +  \rho^2 (1- (\frac{a}{\rho})^2)(1+ (\frac{a}{\rho})^2 + \cdots + (\frac{a}{\rho})^{2(m-2)} + (p+1 -m) (\frac{a}{\rho})^{2(m-1)}) 
                       d \theta^2\\
                      &=&   (1+ \frac{a^2}{\hat r^2 + a^2} + \cdots + (\frac{a^2}{\hat r^2 + a^2})^{(m-2)} + (p+1 -m) (\frac{a^2}{\hat r^2 + a^2})^{(m-1)})^{-1} d\hat r^2\\
                      & & + \frac{\hat r^2}{p^2} (1+ \frac{a^2}{\hat r^2 + a^2} + \cdots + (\frac{a^2}{\hat r^2 + a^2})^{(m-2)} + (p+1 -m) (\frac{a^2}{\hat r^2 + a^2})^{(m-1)}) d \hat \theta^2.
  \eeqn
   When $\hat r \rtarw 0$ , $(1+ \frac{a^2}{\hat r^2 + a^2} + \cdots + (\frac{a^2}{\hat r^2 + a^2})^{(m-2)} + (p+1 -m) (\frac{a^2}{\hat r^2 + a^2})^{(m-1)})\rtarw p.$  
   Therefore, $ds^2|_{\text{fiber}} \rtarw p^{-1} (d\hat r^2 + \hat r^2 d \hat \theta^2).$ The metric restricted to the fiber
 is smooth across the origin. Thus, the metric is complete.

  This process is compatible with the complex structure. 
  We identify the complex plane with the quotient of a fiber through
  $\zeta \rtarw (\frac{c_1}{\sqrt{\sum c_i^2}}, \cdots, \frac{c_m}{\sqrt{\sum c_i^2}}) \zeta^{\frac{1}{p}}.$
  Since  $\frac{y^{m-1}}{y^m + \ga y + \gb} = \frac{1}{p}(y- a^2)^{-1} + l.o.t.,$
  by (\ref{e:u}),  $u = C(y- a^2)^{\frac{1}{p}} + l.o.t. = \hat r^{\frac{2}{p}} + l.o.t.$
  Hence, $|\zeta| = u^{\frac{p}{2}} = \hat r + l.o.t.$  
  Since the Kahler condition is a closed condition, the metric is a complete Kahler metric.

 In special cases when $p = m-1$ or $m$,  $\phi(u)$ can be explicitly  written down.

  $\underline{p= m}.$ The formula (\ref{e:sode}) becomes
     $y^m -a^{2m} = y^{m-1} u y'.$ Therefore, $u = (y^m - a^{2m})^{\frac{1}{m}}$ and hence $y = (a^{2m} + u^m)^{\frac{1}{m}}.$
     We have $\phi'(u) = (1+ a^{2m}u^{-m})^{\frac{1}{m}}.$ The Kahler potential $\phi$ is  the  Calabi's solution \cite{Calabi} :
    $\phi_{c, m} (u) = a^2 (a^{-2m}u^m+ 1)^{\frac{1}{m}} + \frac{a^2}{m} \sum_{j=0}^{m-1} \ln ((a^{-2m} u^m + 1)^{\frac{1}{m}} - \eta^j) \eta^j,$
    where $\eta = e^{2 \pi i /m}$ is the $m$-th unit root. 
    
  $\underline{p= m-1}.$  The formula  (\ref{e:sode}) becomes
     $y^{m-1} -a^{2(m-1)} y  = y^{m-2} u y'.$ Therefore, $u = (y^{m-1} - a^{2(m-1)})^{\frac{1}{m-1}}$ and hence $y = (a^{2(m-1)} + u^{m-1})^{\frac{1}{m-1}}.$
  We have $\phi'(u) = (1+ (a^2 u^{-1})^{m-1})^{\frac{1}{m-1}}.$ The Kahler potential is  $\phi_{c, m-1} (u),$ the Calabi's solution in dimension $m-1$:\\
    $\phi(u) = \phi_{c, m-1}  = a^2 (a^{-2(m-1)}u^{(m-1)}+ 1)^{\frac{1}{m-1}} + \frac{a^2}{m-1} \sum_{j=0}^{m-2} \ln ((a^{-2(m-1)} u^{(m-1)} + 1)^{\frac{1}{m-1}} - \eta^j) \eta^j,$
  where $\eta = e^{2 \pi i /(m-1)}$ is the $(m-1)$-th unit root.
    \epf 
 Finally,  we remark on the volume expansion of the metrics constructed. 
 When $r \rtarw \infty,$
 \beqn
 Vol (\rho \leq r) &=& \int (1+ (m- \frac{3}{2}) \ga \rho^{-2(m-1)}) \rho^{2m-1} d\rho \,d\Omega_{\mathbb{S}^{2m-1}} + l.o.t.\\
 &=& |\mathbb{S}^{2m-1}| (\frac{1}{2m} \rho^{2m}+ \frac{1}{2} (m- \frac{3}{2}) \ga \rho^2 + l.o.t.),
 \eeqn where  $m \geq 2.$
  The coefficient in front of $\rho^2$ changes sign (depending on $p$).
 This is related to the fact that the mass for scalar flat ALE manifolds is not nonnegative  \cite{LeB88}.

\section{Algebraic inequalities for curvature tensors} \label{s:Kato}
 We derive two general algebraic inequalities for Riemannian curvature tensors.
 Then we prove special inequalities for Ricci curvatures of constant scalar curvature  Kahler and (anti-)self-dual 
 metrics.
 
 In this section, the letters $\{i,j, k, l\}$ are indices from $1$ to $n$
    and  $\{a,b,c, d\}$ are from $2$ to $n$ unless otherwise noted.
 
 \begin{lemma} \label{l:wKato} Let $n \geq 3.$ There exists $C=C(n)$ such that for any metric, 
 $$ \sup_{|v|= 1} |\gra_v Rm|^2 \leq \frac{n}{n+2} |\gra Rm|^2 + C |\gd Rm||\gra Rm| + C |\gd Rm|^2.$$
 \end{lemma}
 \bpf
   Let $\{ e_1 =v, \cdots, e_n\}$ be an orthonormal basis at the point.  The curvatures $R_{ijkl, 1}$ consists of three different kinds: 
     $I= \{R_{1a1b,1}, R_{1ab1, 1}, \cdots \}$,  $II= \{R_{1abc,1}, R_{a1bc, 1}, \cdots \}$ and  $III= \{R_{abcd,1}, \cdots \},$  
    where $R_{ijkl}$ has $2, 1$ and $0$  indices of $1$ respectively.
      Note that $$|\gra_1 Rm|^2 = \sum_{a,b,c} 4 |R_{1abc,1}|^2+ \sum_{a,b,c,d}  |R_{abcd,1}|^2+  \sum_{a,b} 4 |R_{1a1b,1}|^2.$$
    
    For $II,$ by
     $
     R_{1abc,1}+ R_{1a1b,c}+ R_{1ac1,b}= 0$ and
    $ R_{1abc,1} + \sum_d R_{dabc,d}= \gd Rm,
 $ we get
   \beqn
    \sum_{a,b,c} 4 |R_{1abc,1}|^2&\leq&  8 \sum_{a,b, c} (|R_{1a1b,c}|^2 + |R_{1ab1,c}|^2) = 4 \sum_{a,b, c} 4 |R_{1a1b,c}|^2, \\
    \sum_{a,b,c} 4 |R_{1abc,1}|^2 & \leq& 4 \sum_{a,b,c} |\sum_d R_{dabc, d}|^2 +  C |\gd Rm||\gra Rm| + C |\gd Rm|^2\\
                & \leq& 4 (n-2) \sum_{a,b,c,d} | R_{dabc, d}|^2 +  C |\gd Rm||\gra Rm| + C |\gd Rm|^2.
   \eeqn   
     Therefore, 
      \beqn
      \sum_{a,b,c} 4 |R_{1abc,1}|^2 &\leq& \frac{3n -8}{4n-8}  4 \sum_{a,b, c} 4 |R_{1a1b,c}|^2+
          \frac{n}{4n-8} 4 (n-2) \sum_{a,b,c,d} | R_{dabc, d}|^2\\
          &+&   C |\gd Rm||\gra Rm| + C |\gd Rm|^2\\
          & \leq& \frac{n}{2}   ( \sum_{a,b, c} 4 |R_{1a1b,c}|^2   +  \sum_{a,b,c,d} 2 | R_{dabc, d}|^2)
          +   C |\gd Rm||\gra Rm| + C |\gd Rm|^2,
        \eeqn where we use that $ \frac{3n -8}{n-2} \leq \frac{n}{2}.$    
       Since $\sum_{a,b,c,d} 2 | R_{dabc, d}|^2 \leq  \sum_{a,b,c,d, e}  | R_{abcd, e}|^2,$ we get
         $$\sum_{a,b,c} 4 |R_{1abc,1}|^2 \leq \frac{n}{2} (\sum_{a,b, c} 4 |R_{1a1b,c}|^2 + 
         \sum_{a,b,c,d, e}  | R_{abcd, e}|^2)+   C |\gd Rm||\gra Rm| + C |\gd Rm|^2.$$

    For $I$ and $III,$  by
     $
     R_{abcd,1}+ R_{ab1c,d}+ R_{abd1,c}= 0$ and
    $R_{1a1b,1} + \sum_c R_{ca1b,c}= \gd Rm,
    $  
    we have
     \beqn
    \sum_{a,b,c,d}  |R_{abcd,1}|^2&\leq& 2 \sum_{a,b,c, d} (|R_{ab1c,d}|^2 + |R_{abd1,c}|^2) =
     \sum_{a,b,c, d} 4 |R_{1abc,d}|^2, \\
    \sum_{a,b} 4 |R_{1a1b,1}|^2 & \leq& 4 \sum_{a,b} |\sum_c R_{ca1b, c}|^2 +  C |\gd Rm||\gra Rm| + C |\gd Rm|^2\\
                & \leq& 4 (n-2) \sum_{a,b,c} | R_{ca1b, c}|^2 +  C |\gd Rm||\gra Rm| + C |\gd Rm|^2.
   \eeqn
     Since $\sum_{a,b,c} | R_{ca1b, c}|^2 \leq \frac{1}{2}  \sum_{a,b,c, d}  |R_{1abc,d}|^2,$ we get
      $$
       \sum_{a,b,c,d}  |R_{abcd,1}|^2+  \sum_{a,b} 4 |R_{1a1b,1}|^2
      \leq   \frac{n}{2} \sum_{a,b,c, d} 4 |R_{1abc,d}|^2+
       C |\gd Rm||\gra Rm| + C |\gd Rm|^2.
     $$   
   
     Combing inequalities of I, II and III and noting that $ \sum_a |\gra_a Rm|^2 = \sum_{a,b, c} 4 |R_{1a1b,c}|^2 + 
         \sum_{a,b,c,d, e}  | R_{abcd, e}|^2+  \sum_{a,b,c, d} 4 |R_{1abc,d}|^2,$ we obtain
     $|\gra_1 Rm|^2 \leq \frac{n}{2} \sum_a |\gra_a Rm|^2 $\\ $+ C |\gd Rm||\gra Rm| + C |\gd Rm|^2.$
      Adding $\frac{n}{2} |\gra_1 Rm|^2$ on both sides  gives the  inequality. 
 \epf
 
 
 \begin{lemma} \label{l:wKato1} Let $n \geq 4.$ There exists $C=C(n)$ such that for any metric, 
 $$\sup_{|v|= 1} |\gra_v W|^2 \leq \frac{n-1}{n+1} |\gra W|^2 + C |\gd W||\gra W| + C |\gd W|^2,$$
 and  therefore,
  $$ \sup_{|v|= 1} |\gra_v Rm|^2 \leq \frac{n-1}{n+1} |\gra Rm|^2 + C |\gra Rc||\gra Rm| + C |\gra Rc|^2.$$
 \end{lemma} 
 \bpf We recall some basic facts about curvatures.  Let $A= \frac{1}{n-2} (Rc - \frac{R}{2(n-1)} g).$
    Then $R_{ijkl} = W_{ijkl} + A_{ik} g_{jl} + A_{jl} g_{ik} - A_{il} g_{jk} - A_{jk} g_{il}.$ In a short hand, 
     we write $Rm = W + A \odot g.$ It is known that $|Rm|^2 = |W|^2 + |A \odot g|^2.$
     Moreover,   $ (\gd W)_{jkl} = \gra^i W_{ijkl}  = (n-3) (A_{jl,k} - A_{jk,l}).$
    Using the second Bianchi, we have
     $$0 = W_{ijkl, m} + W_{ijmk, l} + W_{ijlm, k} + d A \otimes g = W_{ijkl, m} + W_{ijmk, l} + W_{ijlm, k} + \gd W \otimes g,$$
     where $(dA)_{jkl}= A_{jl,k} - A_{jk,l}$ and $d A \otimes g$ represents  tensor products of 
      $dA$ and $g.$
 
   Let $\{ e_1 =v, \cdots, e_n\}$ be an orthonormal basis at the point.     $W_{ijkl, 1}$ consists of three  kinds: 
     $I= \{W_{1a1b,1}, W_{1ab1, 1}, \cdots \}$,  $II= \{W_{1abc,1}, W_{a1bc, 1}, \cdots \}$ and  $III= \{W_{abcd,1}, \cdots \},$  
    where $W_{ijkl}$ has $2, 1$ and $0$  indices of $1$ respectively. Without loss of generality, we may assume
     $W_{1a1b, 1}$ is diagonal. Note that 
     $|\gra_1 W|^2 = \sum_{a,b,c} 4 |W_{1abc,1}|^2+ \sum_{a,b,c,d}  |W_{abcd,1}|^2+  \sum_{a,b} 4 |W_{1a1b,1}|^2.$
    
      In what follows, the summation over $\{ a\neq b \neq c \}$ represents the summation over distinct $a, b, c$ 
         and the summation over $\{a, b, c\}$ represents the summation over all triples $\{a, b,c\}$
       (without the order).  
       
     For II, note that
      $\sum_{a,b,c} 4 |W_{1abc,1}|^2= \sum_{a \neq b \neq c} 4 |W_{1abc, 1}|^2 + \sum_{a \neq b} 8 |W_{1aba, 1}|^2.$
     Since\\
      $  W_{1abc,1}+ W_{1a1b,c}+ W_{1ac1,b}= \gd W \otimes g $ and
     $W_{1abc,1} + \sum_d W_{dabc,d}= \gd W,$
  we get 
   \beqn
   \sum_{a \neq b \neq c} 4 |W_{1abc, 1}|^2 &\leq& \sum_{a\neq b \neq c} 4 |W_{1a1b,c}- W_{1a1c,b}|^2+ C |\gd W||\gra W| + C |\gd W|^2\\
   &=&  
        \sum_{\{a, b, c\}} 8 (|W_{1a1b,c}- W_{1a1c,b}|^2+ |W_{1b1a,c}- W_{1b1c,a}|^2+ |W_{1c1a,b}- W_{1c1b,a}|^2)\\
      & &  + C |\gd W||\gra W| + C |\gd W|^2.
 \eeqn Therefore,
 \beqn
     \sum_{a \neq b \neq c} 4 |W_{1abc, 1}|^2  &\leq & 24 \sum_{\{a, b, c\}}  (|W_{1a1b,c}|^2+  |W_{1a1c,b}|^2+ |W_{1b1c,a}|^2) + C |\gd W||\gra W| + C |\gd W|^2\\
        &=& 3 \sum_{a\neq b \neq c} 4 |W_{1a1b,c}|^2+ C |\gd W||\gra W| + C |\gd W|^2.
   \eeqn
   On the other hand,  for $n \geq 5$
    \beqn \sum_{a \neq b \neq c} 4 |W_{1abc, 1}|^2 &\leq&   \sum_{a\neq b \neq c} 4  |\sum_d W_{dabc,d}|^2 + C |\gd W||\gra W| + C |\gd W|^2 \\
    &\leq& 2(n-2) \sum_{a\neq b \neq c, d} 2  |W_{dabc,d}|^2+ C |\gd W||\gra W| + C |\gd W|^2.
    \eeqn
     Hence, \begin{align} \sum_{a \neq b \neq c} 4 |W_{1abc, 1}|^2 &\leq \frac{3n-7}{4n-8}\, 3 \sum_{a\neq b \neq c} 4 |W_{1a1b,c}|^2+ \frac{n-1}{4n-8}
           2(n-2) \sum_{a\neq b \neq c, d} 2  |W_{dabc,d}|^2 \notag\\
           & + C |\gd W||\gra W| + C |\gd W|^2 \notag\\
      &\leq \frac{n-1}{2}  (\sum_{a \neq b \neq c} 4 |W_{1a1b,c}|^2+  \sum_{a\neq b \neq c, d} 2  |W_{dabc,d}|^2)+  C |\gd W||\gra W| + C |\gd W|^2,  \label{i:kato-a}   
    \end{align} where we use $\frac{3n-7}{4n-8} \cdot 3 \leq \frac{n-1}{2}$ for $n \geq 5.$
    For $n=4,$
    \beqn
    \sum_{a \neq b \neq c} 4 |W_{1abc, 1}|^2 &\leq& 3 \sum_{a\neq b \neq c} 4 |W_{1a1b,c}|^2+ C |\gd W||\gra W| + C |\gd W|^2\\
           &\leq & 3 \sum_{a\neq b \neq c, d} 4 | W_{dadb,c}|^2+ C |\gd W||\gra W| + C |\gd W|^2.
    \eeqn Therefore,
     \beq \sum_{a \neq b \neq c} 4 |W_{1abc, 1}|^2 \leq \frac{3}{2}  (\sum_{a \neq b \neq c} 4 |W_{1a1b,c}|^2+  \sum_{a\neq b \neq c, d} 4  |W_{dadb,c}|^2)+  C |\gd W||\gra W| + C |\gd W|^2.  \label{i:kato-b}   
     \eeq
     
  For $ \sum_{a \neq b} 8 |W_{1aba, 1}|^2,$  by Cauchy inequality we have
    \begin{align}
     \sum_{a \neq b} 8 |W_{1aba, 1}|^2 &\leq  \sum_{a\neq b} 4 |W_{1a1b,a}- W_{1a1a,b}|^2
     +\sum_{a\neq b} 4  |W_{baba,b}+ \sum_{d\neq b}  W_{daba,d}|^2 \notag\\
       &   +  C |\gd W||\gra W| + C |\gd W|^2  \notag\\
    &\leq 4 \sum_{a\neq b}( \frac{n-1}{n-3} (|W_{1a1b,a}|^2+ |\sum_{d\neq b}  W_{daba,d}|^2) +
      \frac{n-1}{2} (|W_{1a1a,b}|^2 +  |W_{baba,b}|^2)) \notag\\
          & +  C |\gd W||\gra W| + C |\gd W|^2. \notag 
  \end{align}
  Thus,
   \begin{align}
     \sum_{a \neq b} 8 |W_{1aba, 1}|^2    &\leq  \frac{n-1}{2} \sum_{a\neq b} (8 |W_{1a1b,a}|^2+  4  |W_{1a1a,b}|^2
            + 4   |W_{baba,b}|^2)+  \frac{n-1}{2} \sum_{d\neq b \neq a} 8  |W_{daba,d}|^2\notag\\
            & +   C |\gd W||\gra W| + C |\gd W|^2. \label{i:kato-c}
    \end{align}
 Combing (\ref{i:kato-a}), (\ref{i:kato-b}) and (\ref{i:kato-c}), we get
    $$\sum_{a,b,c} 4 |W_{1abc,1}|^2 \leq \frac{n-1}{2}(\sum_{a,b,c,d, e}  | W_{abcd, e}|^2
    +  \sum_{a,b,c} 4 |W_{1a1b,c}|^2)+   C |\gd W||\gra W| + C |\gd W|^2.$$

 For I and III, we separate the terms into two types. By
 $$
     W_{abcd,1}+ W_{ab1c,d}+ W_{abd1,c}= \gd W \otimes g, \, \text{and}\;\;
     W_{1a1a,1} + \sum_c W_{ca1a,c}= \gd W,
 $$  we first have that
  $$
   \sum_{\{a,b\} \neq \{c,d\}}  |W_{abcd,1}|^2 \leq 2 \sum_{\{a,b\} \neq \{c,d\}}  (|W_{ab1c,d}|^2 + |W_{abd1,c}|^2)
   = \sum_{\{a,b\} \neq \{c,d\}} 4 |W_{1cab,d}|^2.   
   $$
 The remaining terms in III plus the terms in I are
 \beqn
      & &\sum_{a} 4 |W_{1a1a,1}|^2 +  \sum_{\{a,b\} = \{c,d\}}  |W_{abcd,1}|^2 \\
      & \leq& 4 \sum_a |\sum_c W_{ca1a, c}|^2 +  4 \sum_{a < b} |W_{abab,1}|^2
      +C |\gd W||\gra W| + C |\gd W|^2\\
                & \leq&  4 \sum_a |\sum_c W_{ca1a, c}|^2 +  4 \sum_{a < b} |W_{ba1a,b}+ W_{ab1b,a}|^2
      +C |\gd W||\gra W| + C |\gd W|^2.
  \eeqn
 Define an $(n-1) \times (n-1)$ matrix by $M_{ab}= W_{ab1b,a}.$  
 Then $M_{ab}= 0$ when $a=b$ and $\sum_b M_{ab}= 0$ for all $a.$
 
 Claim: $\sum_a |\sum_c M_{ca}|^2 +   \sum_{a < b} |M_{ba}+ M_{ab}|^2 \leq (n-1) \sum_{a,b} |M_{ab}|^2.$\\
 Assuming the claim, we get that 
 $$
  4 \sum_a |\sum_c W_{ca1a, c}|^2 +  4 \sum_{a < b} |W_{ba1a,b}+ W_{ab1b,a}|^2
  \leq  \frac{n-1}{2}  \sum_{a,b} 8 |W_{1bab,a}|^2.
 $$  Note  that $\sum_a |\gra_a W|^2 = \sum_{a,b, c} 4 |W_{1a1b,c}|^2 + 
         \sum_{a,b,c,d, e}  | W_{abcd, e}|^2+  \sum_{a,b,c, d} 4 |W_{1abc,d}|^2$ and 
         $\sum_{\{a,b\} \neq \{c,d\}} 4 |W_{1cab,d}|^2+  \sum_{a,b} 8 |W_{1bab,a}|^2 \leq \sum_{a,b,c, d} 4 |W_{1abc,d}|^2.$
  Now combing inequalities of I, II and III, 
         we obtain
     $|\gra_1 W|^2 \leq \frac{n-1}{2} \sum_a |\gra_a W|^2 + C |\gd W||\gra W| + C |\gd W|^2,$
     which proves the first part of the lemma. The second part follows by the curvature decomposition
     $|\gra Rm|^2 = |\gra W|^2 + |\gra A \odot g|^2.$

 It remains to prove the Claim. We use Lagrange multiplier to estimate the maximum of
    $F= \sum_a |\sum_c M_{ca}|^2 +   \sum_{a < b} |M_{ba}+ M_{ab}|^2$ under the constraints
    $H_0 =  \sum_{a,b} |M_{ab}|^2-1 = 0$ and $H_a= \sum_b M_{ab}= 0$ for $(n-1) \times (n-1)$ matrices $M_{ab}$
    satisfying $M_{ab}= 0$ when $a=b.$ At a critical point, we have $\gra F = \gl_0 \gra H_0 + \sum_a \gl_a \gra H_a.$ 
    Therefore,
    \beq \label{e:lag}
     M_{ab}+ M_{ba} + \sum_c M_{cb}= \gl_0 M_{ab} + \frac{\gl_a}{2}.
   \eeq
    Fixing $a$ and using $H_c= 0,$ we get
    $ \sum_b M_{ab}+ \sum_bM_{ba} + \sum_{b \neq a} \sum_c M_{cb}
    = \sum_b M_{ba}- \sum_c M_{c a}=    (n-1) \frac{\gl_a}{2}.$
   Hence, $\gl_a = 0$ for all $a.$ Thus, $\gra F = \gl_0 \gra H_0.$ This is the equation for the quadratic polynomial  
   $F$ under the constraint $H_0 = 0.$ Therefore, the maximum of $F$ is the maximum of $\gl_0.$
   Going back to (\ref{e:lag}), we have
   $ \sum_a M_{ab}+ \sum_a M_{ba} + \sum_{a \neq b} \sum_c M_{cb}= (n-1) \sum_a M_{ab} = \gl_0 \sum_a M_{ab}.$
   Hence, either $\gl_0 = n-1$ or $\sum_a M_{ab}= 0$ for all $b.$
  If $\sum_a M_{ab}= 0$ for all $b,$ then $F = \sum_{a < b} |M_{ba}+ M_{ab}|^2 \leq 2\sum_{a,b} |M_{ab}|^2 = 2< n-1.$ 
  Thus, we conclude that $F \leq n-1.$
 \epf
 We remark that  for Ricci flat manifolds, it was shown in \cite{BKN89}, \cite{Bran00} that $|\gra |W||^2 \leq \frac{n-1}{n+1} |\gra W|^2.$

  Now we study special cases when $g$ is Kahler or self-dual with constant $R.$  
  We will show that  $\sup_{|v|= 1}  |\gra_v Rc|^2 \leq \frac{n}{n+2} |\gra Rc|^2,$ 
  which is stronger than   $|\gra |Rc||^2 \leq \frac{n}{n+2} |\gra Rc|^2.$
  It was shown in \cite{TV05} that in dimension 4  for metrics satisfying $\gd W^+=0$ with constant $R,$
   it holds   $|\gra |Rc||^2 \leq \frac{2}{3} |\gra Rc|^2.$
  
   \begin{lemma} \label{l:kahler} Let $n=2m \geq 4.$ Suppose $g$ is Kahler and $R$ is constant. Then
   $$\sup_{|v|= 1}  |\gra_v Rc|^2 \leq \frac{n}{n+2} |\gra Rc|^2. $$
 \end{lemma}
  \bpf 
  Let $\{ e_1 =v, \cdots, e_m, J e_1, \cdots, J e_m\}$ be an orthonormal basis at the point. The letters $i,j, k, l$ are indices from $1$ to $m$
    and  $a,b,c, d$ are from $2$ to $m.$  The curvatures $\gra_1 Rc$ consists of three different kinds (if $m \geq 3$  only the first two kinds):
     $I= \{R_{ii,1} \},$  $II = \{ R_{1 a, 1}, R_{1 Ja, 1}\},$ and  $III = \{ R_{a b, 1}, R_{a Jb, 1}, \text{for} \, a<b \}.$  
         Note that  $|\gra Rc|^2 = \sum_i |\gra_i Rc|^2+ \sum_{Ji} |\gra_{Ji} Rc|^2$ and 
     $$|\gra_1 Rc|^2 = \sum_{ij} 2 (|R_{ij,1}|^2 +  |R_{i Jj, 1}|^2) =  \sum_i 2 |R_{ii, 1}|^2   
         +  \sum_{1 \leq i < b} 4 (|R_{i b, 1}|^2+ |R_{i Jb, 1}|^2).$$
     
   For I, by (\ref{e:K-bianchi})
  and  $R_{Z_a \bar Z_a, Z_1} = R_{Z_1 \bar Z_a, Z_a},$ we get 
   $$
     R_{aa,1}= R_{1a,a}+ R_{1 Ja,Ja}, \, \text{and}\;\;
     R_{aa,J1} =  R_{1 a,Ja}-  R_{1 Ja. a}.
   $$ 
    Therefore, using $\sum_i R_{ii}$ is constant we have
    \beqn
     \sum_i 2 |R_{ii,1}|^2&=& 2 |\sum_a R_{aa,1}|^2 + 2\sum_a |R_{aa,1}|^2 \leq  2 m\sum_a |R_{aa,1}|^2\\
      &=& 2m \sum_a |R_{1a,a}+ R_{1 Ja,Ja}|^2 \leq m \sum_a 4 (|R_{1a, a}|^2 + |R_{1 Ja, Ja}|^2).
 \eeqn
    
   For II and III, let $1 \leq i < b \leq m.$
  By  (\ref{e:K-bianchi}) and $R_{Z_b \bar Z_i, Z_1} = R_{Z_1 \bar Z_i, Z_b},$ we get
    $$
     R_{i b,1}= R_{i Jb, J1}+ R_{1 i,b} + R_{1 Ji, Jb}, \, \text{and}\;\;
     R_{i Jb ,1} = - R_{i b,J1}-  R_{1 Ji, b} + R_{1 i, Jb}.
   $$ 
 Therefore,
  $$
     R_{1 b,1}= R_{1 Jb, J1}+ R_{1 1,b}, \, \text{and}\;\;
     R_{1 Jb ,1} = - R_{1 b,J1} + R_{1 1, Jb}.
  $$ Using  $|R_{11, b}|^2 \leq (m-1) \sum_a |R_{aa,b}|^2$ and  $8\frac{m-1}{m} \leq 2m,$ II becomes
   \beqn
     & &\sum_a 4 (|R_{1a,1}|^2+ |R_{1 Ja, 1}|^2)  
      \leq 8 \sum_a ( |R_{11,a}|^2+ |R_{1 1,Ja}|^2 + |R_{1 Ja, J1}|^2 + |R_{1 a, J1}|^2)\\
      &\leq& 8 \frac{m-1}{m} \sum_{i, a} (|R_{ii,a}|^2+ |R_{ii,Ja}|^2) + 8 \sum_a (|R_{1 Ja, J1}|^2 + |R_{1 a, J1}|^2)\\
      &\leq& m \sum_{i, a} 2(|R_{ii,a}|^2+ |R_{ii,Ja}|^2) + m \sum_a 4(|R_{1 Ja, J1}|^2 + |R_{1 a, J1}|^2).
 \eeqn 
  Similarly, for $1<a < b \leq m,$
    $$
     R_{a b,1}= R_{a Jb, J1}+ R_{1 a,b} + R_{1 Ja, Jb},\, \text{and}\;\;
     R_{a Jb ,1} = - R_{a b,J1}-  R_{1 Ja, b} + R_{1 a, Jb}.
    $$  And III becomes
     \beqn
     \sum_{a< b} 4(|R_{a b,1}|^2+ |R_{a Jb, 1}|^2) =  \sum_{a< b} 4( |R_{a Jb,J1}+ R_{1a, b}+ R_{1 Ja, Jb}|^2 + |R_{a b,J1}+ R_{1 Ja, b}- R_{1a, Jb}|^2) \\
      \leq m \sum_{a<b} 4( |R_{a Jb,J1}|^2+ |R_{a b,J1}|^2 + |R_{1 a, b}|^2 + |R_{1 Ja, b}|^2+ |R_{1 Ja, Jb}|^2 + |R_{1 a, Jb}|^2),
     \eeqn where we use $m \geq 3.$ (Note that when $m =2,$ there is no terms in III.)

 Combining I, II and III (when $m=2,$ combining only I and II), 
   \beqn
    |\gra_1 Rc|^2
     &\leq& m ( \sum_{i, a} 2(|R_{ii,a}|^2+ |R_{ii,Ja}|^2)
      +  \sum_{a, b} 4 (|R_{1a, b}|^2 + |R_{1 Ja, b}|^2)
     \\
     &+& \sum_{1 \leq i < b, j} 4(|R_{i b, Jj}|^2+ |R_{i Jb, Jj}|^2))
     \leq  m (\sum_a |\gra_a Rc|^2+ \sum_i |\gra_{Ji} Rc|^2).
  \eeqn 
 Adding $m |\gra_1 Rc|^2$ on both sides  gives the  inequality. 
  \epf
 
 We consider a weaker condition $\gd W^-=0$ ($\gd W^+=0$) than (anti-)self-dual.
  \begin{lemma} \label{l:self-dual}
 Let $n=4.$
  Suppose $\gd W^+=0$ (or $\gd W^-=0$),  and $R$ is constant. Then
 $$\sup_{|v|= 1} |\gra_v Rc|^2 \leq \frac{2}{3} |\gra Rc|^2. $$
 \end{lemma}
 \bpf 
 We will prove the case when $\gd W^+=0.$  The proof for the case $\gd W^-=0$ is similar.
 Let $\{ e_1 =v, e_2, e_3, e_4\}$ be an
 orthonormal basis at the point. The  indices $i,j, k$ are from $1$ to $4$
    and  $a,b,c$ are from $2$ to $4.$ 
 Without loss of generality, we may assume
  $R_{ab, 1}$ is diagonal. Hence, $R_{ij, 1}$  consists of two
 different kinds: $I= \{R_{ii,1}\}$ and
 $II= \{R_{1a,1}, R_{a1,1}\}.$

 Since $\gd W^+ = 0,$ we have
 $\frac{1}{2} d Rc = d A = \gd W = \gd W^- \in T^*M \otimes \Lambda^2_-.$
 Therefore, $d Rc$ is perpendicular to $\Lambda^2_+.$ Recall that
 a basis for $\Lambda^2_+$ is $\{\frac{1}{\sqrt 2} (e_1 \wedge e_2+ e_3 \wedge e_4),
  \frac{1}{\sqrt 2}(e_1 \wedge e_3- e_2 \wedge e_4),
  \frac{1}{\sqrt 2}(e_1 \wedge e_4+ e_2 \wedge e_3)\}.$
 Hence, we obtain the following twelve equations:
 \begin{equation} \label{e} \left\{
 \begin{array}{l}
 R_{i2, 1} = R_{i1, 2} + R_{i3, 4} -R_{i4, 3},\\
 R_{i3, 1} = R_{i1, 3} + R_{i4, 2} - R_{i2, 4},\\
 R_{i4, 1} = R_{i1, 4} + R_{i2, 3} - R_{i3, 2},
 \end{array}\right .
 \end{equation}
  for $i= 1, \cdots, 4.$

 For terms in I, by (\ref{e}) we get
  $$\begin{array} {ll}
 R_{22, 1} &= R_{12, 2} + R_{23, 4} -R_{24, 3},\\
 R_{33, 1} &= R_{13, 3} + R_{34, 2} - R_{32, 4},\\
 R_{44, 1} &= R_{14, 4} + R_{42, 3} - R_{43, 2},\\
 R_{11, 1} &= -R_{22, 1} - R_{33, 1} - R_{44, 1}
           = -R_{12, 2} - R_{13, 3} - R_{14, 4}.
 \end{array}$$
Therefore,
 \begin{eqnarray*}
  && |R_{11, 1}|^2+ |R_{22, 1}|^2 + |R_{33, 1}|^2 + |R_{44,
  1}|^2\\
   &=& 4 (|R_{12, 2}|^2 + |R_{13, 3}|^2 + |R_{14, 4}|^2 +|R_{23, 4}|^2 +|R_{24,
   3}|^2 + |R_{34,2}|^2)- (R_{23, 4} +R_{24,3} + R_{34,2})^2\\
   && - (R_{12, 2} - R_{13, 3} - R_{23,
  4})^2- (R_{13, 3} - R_{14, 4} - R_{34, 2})^2 -(R_{14, 4} - R_{12, 2} - R_{24,
  3})^2\\
  &\leq& 4 (|R_{12, 2}|^2 + |R_{13, 3}|^2 + |R_{14, 4}|^2 +|R_{23, 4}|^2 +|R_{24,
   3}|^2 + |R_{34,2}|^2).
 \end{eqnarray*}

   For II, using (\ref{e}) again and the Bianchi identity
  $R_{21,1} =  - R_{22,2}- R_{23,3}- R_{24,4},$ we get
  \begin{eqnarray*}
 & & |R_{12, 1}|^2 + |R_{21, 1}|^2 =
  |R_{11, 2}+ R_{13,4}- R_{14,3}|^2 + |R_{22,2}+ R_{23,3}+
  R_{24,4}|^2\\
  &=& 2(|R_{11,2}|^2+ |R_{22,2}|^2) + 4(|R_{13,4}|^2+ |R_{14,3}|^2+|R_{23,3}|^2 +
  |R_{24,4}|^2) - 2(R_{13, 4}+ R_{14, 3})^2 \\
  & &- 2(R_{23,3}- R_{24,2})^2- (R_{11, 2}+ R_{14,3}-  R_{13,4})^2 - (R_{22, 2}- R_{23,3}-
  R_{24,4})^2
 \\
  &\leq& 2(|R_{11,2}|^2+ |R_{22,2}|^2) + 4(|R_{13,4}|^2+ |R_{14,3}|^2+|R_{23,3}|^2 +
  |R_{24,4}|^2).
  \end{eqnarray*}
 Similarly, we also have
 $$\begin{array}{l}
   |R_{13, 1}|^2 + |R_{31, 1}|^2 \leq 2(|R_{11,3}|^2+ |R_{33,3}|^2) +
  4(|R_{12,4}|^2+ |R_{14,2}|^2+|R_{32,2}|^2 +|R_{34,4}|^2),\\
  |R_{14, 1}|^2 + |R_{41, 1}|^2 \leq 2(|R_{11,4}|^2+ |R_{44,4}|^2) +
  4(|R_{12,3}|^2+ |R_{13,2}|^2+|R_{42,2}|^2 +|R_{43,3}|^2).
  \end{array}$$
 Adding the above three inequalities and the one from I together, we 
 get $|\gra_1 Rc|^2 \leq 2 \sum_a |\gra_a Rc|^2.$
 Therefore, $|\gra_1 Rc|^2 \leq \frac{2}{3}  |\gra Rc|^2.$
 \epf


\section{Proofs of Theorem~\ref{t:delRm}, \ref{t:Kahler} and Proposition~\ref{p:poh}}
 \label{s:thm}

\bpf[Proof of Proposition~\ref{p:poh}]
 Note that $\gra_i X^j = \gd_{ij}+ \Gamma \ast X.$ We compute
    \beqn
  & &\int_{\gO} \li \Delta T, X^i \gra_i T \ri dV
   =- \int_{\gO}  \li \gra_j T,  \gra_j (X^i \gra_i T) \ri dV + \int_{\de \gO}  \li  \gra_{\nu} T,   X^i \gra_i T \ri d\gs\\
   &=&- \int_{\gO}  \li \gra_j T,  \gra_j T + \gC \ast X \ast \gra T + X^i \gra_j \gra_i T \ri dV
     +  \int_{\de \gO}  \li  \gra_{\nu} T,   X^i \gra_i T \ri d\gs, 
   \eeqn where $\nu$ is the unit outer normal.
   Exchanging the derivatives, we get
   \beqn
   \int_{\gO} \li \Delta T, X^i \gra_i T \ri dV
   &\leq& - \int_{\gO}  \li \gra_j T,  \gra_j T + \gC \ast X \ast \gra T + X^i \gra_i \gra_j T + X \ast Rm \ast T \ri dV\\
      &+&  \int_{\de \gO}  \li  \gra_{\nu} T,   X^i \gra_i T \ri d\gs\\
      &=&  \int_{\gO} (- |\gra T|^2 - \frac{1}{2}  X^i \gra_i  |T|^2 + \gra T \ast \gra T\ast \gC \ast X \\
      &+& \gra T\ast T \ast X \ast Rm ) dV +  \int_{\de \gO}  \li  \gra_{\nu} T,   X^i \gra_i T \ri d\gs.
    \eeqn
  Therefore,
   \beqn
  \int_{\gO} \li \Delta T, X^i \gra_i T \ri dV 
  &=& \int_{\gO}  \left( \frac{n-2}{2}  |\gra T|^2 
  +   \gra T \ast \gra T \ast \gC \ast X +   \gra T \ast  T \ast Rm \ast X \right)  dV\\
  &+&  \int_{\de \gO}  \li  \gra_{\nu} T,   X^i \gra_i T \ri d\gs- \int_{\de \gO} \frac{1}{2}  \li X, \nu \ri |\gra T|^2 d\gs.
 \eeqn
\epf

  Now we are ready to prove the main theorems.
 \bpf[Proof of Theorem~\ref{t:delRm}]
      By Bando-Kasue-Nakajima \cite{BKN89}, if $|Rm| = O(r^{-(2+ \ga)}),$ then there exist coordinates of
       order $\ga.$ The problem boils down to estimating the decay of $|Rm|.$
     
       Let $D_r$ be  the complement of the geodesic ball of radius $r.$
  Applying Lemma~\ref{l:inhomo} to (\ref{e:curv}) we have $|Rm| \in L^{\frac{n}{2} \gc}.$ 
 Moreover, by Lemma~\ref{l:inhomo} (a) with $p = \frac{2n}{n-2}$ and    the Sobolev inequality, 
 \beq\label{i:Rm}
 \sup_{D_{2r}} |Rm| \leq Cr^{-\frac{n-2}{2}} \|Rm\|_{L^\frac{2n}{n-2}(D_r)} + C r^{-(s+2)} \leq Cr^{-\frac{n-2}{2}} \|\gra Rm\|_{L^2(D_r)} + C r^{-(s+2)}. 
 \eeq
  The formula (\ref{i:Rm})   passes the decay of $ \int_{D_r}  |\gra Rm|^2   dV$ to that of $|Rm|.$
   We will use Pohozaev's identity and apply ODE estimates (Lemma~\ref{l:ode}) to 
  $f(r) =  \int_{D_r}  |\gra Rm|^2  dV$  to get  the decay of $f(r).$ Then apply the above regularity.

  (a)   Case 1.  $n \geq 5.$
   
   By Lemma~\ref{l:inhomo}, we get that $|Rm|= O(r^{-\ga_0})$ for all $\ga_0 < \min \{n-2, s+2\}.$ Therefore, $\ga_0 > 0.$
   By Proposition~\ref{t:reg-curv}, we have $|\gra^k Rm|= O(r^{-(\ga_0 +k)}).$ By \cite{BKN89}, 
     there exist coordinates $x$ and 
    $0<\gd_0<1$ such that  $C^{-1} r \leq |x| \leq C r,$ $|g-\gd| = O(|x|^{-\gd_0})= O(r^{-\gd_0})$ and $|\de g| = O(|x|^{-\gd_0 -1})=
    O(r^{-\gd_0 -1}),$
    where $\gd$ is the standard metric on the Euclidean space.  
     Since distances $|x|$ and $r$ are equivalent, we will still denote $|x|$ by $r$ and denote
     by $D_r$ the complement of the coordinates ball of radius $r$ and by $S_r$
     the coordinates sphere of radius $r.$

       Applying Proposition~\ref{p:poh} with
     $T= Rm$ and $\Omega= D_r$ and using (\ref{e:curv}), we obtain
    \beqn 
    & &\int_{D_r} (\gra \gd Rm \ast X \ast \gra Rm + Rm \ast Rm \ast X \ast \gra Rm) dV= \int_{D_r}  \Delta R_{ijkl} X_p R_{ijkl,p}  dV\\
    &\geq& \int_{D_r} ( \frac{n-2}{2} |\gra Rm|^2 - C|\gra Rm|^2 |\gC|  |X| -  C|\gra Rm|   |Rm|^2  |X| ) dV\\
  &+& \int_{S_r}   \nu_p R_{ijkl, p} X_q R_{ijkl,q} d\gs-\int_{S_r} \frac{1}{2}  X_k \nu_k  |\gra Rm|^2 d\gs,
 \eeqn where $\nu$ is the unit outward normal.  
   Note that  $\nu= - \frac{X}{r} + O(r^{-\gd_0}).$ Hence,
   \begin{align}
 \int_{D_r}  \frac{n-2}{2} |\gra Rm|^2   dV&\leq    \int_{D_r} C(|\gra Rm|^2 |\gC|  |X| 
    +  |\gra Rm|   |Rm|^2  |X| + |\gra \gd Rm| |X| |\gra Rm|) dV \notag \\ 
   &+   r \int_{S_r}   |\gra_{\nu} Rm|^2 d\gs    -r \int_{S_r} \frac{1}{2}  |\gra Rm|^2 d\gs + C \int_{S_r} |\gra Rm|^2 r^{1-\gd_0} d\gs.    \label{i:poh-main}
   \end{align}
 By Lemma~\ref{l:wKato}, the above formula becomes
    \begin{align}
 & \int_{D_r}  \frac{n-2}{2} |\gra Rm|^2   dV \leq       \int_{D_r} C(|\gra Rm|^2 |\gC|  |X| 
    +  |\gra Rm|   |Rm|^2  |X| + |\gra \gd Rm| |X| |\gra Rm|) dV     \notag  \\
  & + \left( \frac{n}{n+2}- \frac{1}{2} + C r^{- \gd_0} \right) r \int_{S_r}   |\gra Rm|^2 d\gs +   r \int_{S_r} C  (|\gd Rm|^2 + |\gd Rm||\gra Rm|) d\gs. \label{i:poh-1}
   \end{align}
  Using  $|\gra^k \gd Rm| = O(r^{-(k+s+3)})$ and  $|\gra Rm| = O(r^{-(\ga_0 + 1)}),$ 
  we have
 \beqn
 \int_{D_r} \left( \frac{n-2}{2} - \epsilon \right) |\gra Rm|^2   dV 
 &\leq&    \left(\frac{n-2}{2(n+2)} + \epsilon \right) r \int_{S_r}   |\gra Rm|^2 d\gs\\
 &+&      O( r^{-(2s +6 -n)}+r^{-(4 \ga_0-2-n)}+ r^{-(2 \ga_0 + 2 -n + \gd_0)}).
   \eeqn 
  Hence, for some $\epsilon'$ small
   \beq\int_{D_r}  |\gra Rm|^2   dV 
 \leq (n+2- \epsilon')^{-1}   r \int_{S_r}   |\gra Rm|^2 d\gs 
     +      O( r^{-(2s +6 -n)}+r^{-(4 \ga_0-2-n)}+  r^{-(2 \ga_0 + 2 -n + \gd_0)}). \label{i:ode1}
    \eeq 
      
  When $s \geq n-2,$ we have    $2 \ga_0 + 2 -n + \gd_0> 2s+6-n\geq n+2.$ If $n \geq 6,$ then $4 \ga_0 -2 -n = 3n-10 -\tilde \gep > n+2 - \gep'$
  by choosing $\tilde \gep \ll \gep'.$ Hence, by Lemma~\ref{l:ode} with $f(r) = \int_{D_r}  |\gra Rm|^2   dV,$ 
  we have $\int_{D_r}  |\gra Rm|^2   dV= O(r^{-(n+2 - \gep')}).$ 
  Since  $|x|$ and the geodesic distance are equivalent,
   we can apply (\ref{i:Rm}) to the complement of the coordinates ball $D_r.$
  Thus,  $|Rm| =O( r^{-(n-\frac{\gep'}{2})}+ r^{-(s+2)})=O( r^{-(n-\frac{\gep'}{2})}).$ Therefore, $|Rm| = O(r^{-\ga})$ for all $\ga < n.$
 If $n=5,$ $4 \ga_0 -2 -n = 5 -\tilde \gep < n+2 - \gep'= 7- \gep'.$ Again, by Lemma~\ref{l:ode} we get
  $\int_{D_r}  |\gra Rm|^2   dV= O(r^{-(5 - \tilde \gep)}).$ By   (\ref{i:Rm}), 
  $|Rm| =O( r^{-(4-\frac{\tilde \gep}{2})}).$   Go back to (\ref{i:ode1}) now with $\ga_0 =4-\frac{\tilde \gep}{2}.$
   Then $4 \ga_0 -2 -n = 9 -2\tilde \gep > n+2 - \gep'= 7- \gep'.$ By  Lemma~\ref{l:ode} and (\ref{i:Rm}), we finally obtain
   $|Rm|= O(r^{-(5 - \frac{\gep'}{2})})$ and  hence $|Rm| = O(r^{-\ga})$ for all $\ga < 5.$
    This completes  the proof for $s=n-2.$
      
    When $s> n-2,$ by (\ref{i:poh-1}) with $|Rm| = O(r^{-\ga})$ and  $|\gra Rm|= O(r^{-(\ga + 1)})$ (by Proposition~\ref{t:reg-curv}),
    we have 
    \beqn
     \int_{D_r}  \frac{n-2}{2} |\gra Rm|^2   dV &\leq&       
   \left( \frac{n}{n+2}- \frac{1}{2} \right) r \int_{S_r}   |\gra Rm|^2 d\gs\\
   & +& O(r^{-(2 \ga + 2 + \gd_0 -n)}+ r^{-(2s+6-n)}+r^{-(\ga + s+ 4-n)}
   + r^{-(3 \ga -n)}).
    \eeqn  
  Let $\ga = n- \tilde \gep$ be close to $n.$ Then $2 \ga + 2 + \gd_0 -n, 2s+6-n, \ga + s+ 4-n$ and  $3 \ga -n$
     are all strictly larger than $n+2.$ Hence, $ \int_{D_r}  |\gra Rm|^2   dV \leq       
   (n+2)^{-1} r \int_{S_r}   |\gra Rm|^2 d\gs+ O(r^{- \gb})$ for some $\gb > n+2.$ By Lemma~\ref{l:ode} and (\ref{i:Rm}) again,
    $|Rm| = O(r^{-n}).$
    
  When $n-4\leq s < n-2,$ we have $2 \ga_0 + 2 -n + \gd_0,  n+2 - \gep' > 2s+6 -n.$ 
  If $n \geq 6,$ then $4 \ga_0 -2 -n = 3n -10 -\tilde \gep > 2s+6-n.$    
   Applying     Lemma~\ref{l:ode}  to (\ref{i:ode1}) gives $\int_{D_r}  |\gra Rm|^2   dV= O(r^{-(2 s + 6-n)}).$ Thus, 
   $|Rm| = O(r^{-(s+2)}).$ If $n=5$ and $2s+1 < 5,$ then $2 s + 6-n< 4 \ga_0 -2 -n = 5 -\tilde \gep.$ 
    By  Lemma~\ref{l:ode}, (\ref{i:ode1}) and  (\ref{i:Rm}), we get $|Rm| = O(r^{-(s+2)}).$   
    If $n=5$ and $2s+1 \geq  5,$ then $2 s + 6-n> 4 \ga_0 -2 -n = 5 -\tilde \gep.$  
      Hence, $|Rm| = O(r^{-(4 -  \frac{\tilde \gep}{2})}).$ Applying the same argument to 
      (\ref{i:ode1}) with $\ga_0 = 4 -  \frac{\tilde \gep}{2},$ we finally arrive at $|Rm| = O(r^{-(s+2)}).$ 
      
   When $s < n-4,$ we have $ n+2 - \gep', 2 \ga_0 + 2 -n + \gd_0> 2s+6 -n$ and $4 \ga_0 -2 -n = 4s+6-n -\tilde \gep > 2s+6-n.$   
    By  Lemma~\ref{l:ode}, (\ref{i:ode1}) and  (\ref{i:Rm}), we get $|Rm| = O(r^{-(s+2)}).$ 
           
 
 Case 2. $n=4.$
 
  Consider the equation (\ref{e:curv}).
  By Lemma~\ref{l:inhomo} (a) and Proposition~\ref{t:reg-curv}, we have
  \begin{align}
  |Rm| \leq Cr^{-4/p} \|Rm\|_{L^p} + Cr^{-(2+s)},\notag \\
 \int_{B_{r}} |\gra Rm|^2 dV \leq C r^{-2} \int_{B_{2r}} |Rm|^2 dV + C r^{-2(s+1)}, \label{i:dim4}
 \end{align} for $p> 1.$
  Since $|Rm|^2$ is integrable, let $\int_{D_r} |Rm|^2 dV= \gep_0(r)^2,$ where $\gep_0 (r) \rtarw 0$ when $r \rtarw \infty.$
  Therefore, by \cite{TV06}, 
    we have $|g- \gd| \leq C \gep_0 (r), |\de g| \leq C \gep_0 (r)/r$ and $ |Rm| \leq C \gep_0 (r)/r^2.$
 
 By (\ref{i:poh-1}), (noting that $\nu = -\frac{X}{r}+ o(1),$ then $r^{-\gd_0}$ becomes $\epsilon_0 (r)$ in (\ref{i:poh-1}))
    \begin{align} 
& \int_{D_r}   |\gra Rm|^2   dV \leq       \int_{D_r} C(|\gra Rm|^2 |\gC|  |X| 
    +  |\gra Rm|   |Rm|^2  |X| + |\gra \gd Rm| |X| |\gra Rm|) dV    \notag   \\
  & + (\frac{1}{6}+ C\epsilon_0 (r)) r \int_{S_r}   |\gra Rm|^2 d\gs +   r \int_{S_r} C  (|\gd Rm|^2 + |\gd Rm||\gra Rm|) d\gs. \label{i:poh-dim4} 
   \end{align}
  Let $A_r= \{x: r/2 < |x|< 3r/2\}.$ By (\ref{i:dim4}) and Sobolev inequality,
  \beqn
  \int_{A_r}  |\gra Rm|   |Rm|^2 r dV &\leq& C r \|\gra Rm\|_{L^2(A_r)} \|Rm\|_{L^4(A_r)}^2 \leq C r \|\gra Rm\|_{L^2(A_r)}^3\\
   &\leq & C( \|Rm\|_{L^2(A_r)} + r^{-s})  \|\gra Rm\|_{L^2(A_r)}^2. 
  \eeqn
  Hence,  $$ \int_{D_r}  |\gra Rm|   |Rm|^2 |X| dV \leq  C(  \gep_0 (r) + r^{-s})  \|\gra Rm\|_{L^2(D_r)}^2.$$
  Now (\ref{i:poh-dim4})  becomes
  \beqn
  (1-\gep_0 (r)- \gep)  \int_{D_r}   |\gra Rm|^2  dV & \leq&  C \int_{D_r} |\gra \gd Rm|^2 |X|^2 dV\\
   &+& (\gep+ \frac{1}{6}+ C\epsilon_0 (r)) r \int_{S_r}   |\gra Rm|^2 d\gs 
  +   r \int_{S_r} C  |\gd Rm|^2 d\gs 
  \eeqn for some $\gep > 0$ small.
  Thus, by $|\gra^k \gd Rm| = O(r^{-(s+3+k)})$ for some $\epsilon'$ small
   \beq\int_{D_r}  |\gra Rm|^2   dV 
 \leq (6- \epsilon')^{-1}   r \int_{S_r}   |\gra Rm|^2 d\gs 
     +      O( r^{-(2s +2)}). \label{i:ode2}
    \eeq 
  
   When $s\geq 2,$ we have    $2s+2\geq 6.$
  Hence, by Lemma~\ref{l:ode}, 
  we have $\int_{D_r}  |\gra Rm|^2   dV= O(r^{-(6 - \gep')}).$ By (\ref{i:Rm}), 
  $|Rm| =O( r^{-(4-\frac{\gep'}{2})}).$ Therefore, $|Rm| = O(r^{-\ga})$ for all $\ga < 4.$
  This completes the proof for $s= 2$ case.
   When $s> 2,$ now  $|Rm| = O(r^{-\ga})$ for $\ga= 4- \tilde \gep.$
    Thus, by \cite{BKN89} there exist coordinates such that  
   $|g-\gd| = O(r^{-\gd_0})$ and $|\de g| = O(r^{-\gd_0 -1})$ for some fixed $\gd_0.$
    The rest of the argument is exactly the same as in Case 1 ($n \geq 5$),  $s > n-2.$ 
    We have
    $$
     \int_{D_r}   |\gra Rm|^2   dV \leq       
    \frac{1}{6} r \int_{S_r}   |\gra Rm|^2 d\gs
    + O(r^{-(2 \ga - 2 + \gd_0)}+ r^{-(2s+2)}+r^{-(\ga + s)}
   + r^{-(3 \ga -4)}).
    $$ 
 We can check that $2 \ga - 2 + \gd_0, 2s+2, \ga + s$ and  $3 \ga -4$ are 
     all strictly larger than $6.$ Hence,  By Lemma~\ref{l:ode} and (\ref{i:Rm}),
    $|Rm| = O(r^{-4}).$
 
  When $s < 2,$ we have $2s+2 < 6 - \gep'.$   
   Applying     Lemma~\ref{l:ode}  to (\ref{i:ode2}) gives $\int_{D_r}  |\gra Rm|^2   dV= O(r^{-(2 s + 2)}).$ Thus, 
   $|Rm| = O(r^{-(s+2)}).$ 
 
 (b) Since $\gd Rm = d Rc,$ $|\gra^k \gd Rm| \leq |\gra^{k+1} Rc|.$ And  by (a) we already know that $|g-\gd| = O(r^{-\gd_0}),$ $|\de g| = O(r^{-\gd_0 -1})$
    and $|Rm| = O(r^{-\ga})$ for all $\ga < n.$ We will prove that under assumptions $|\gra Rc| = O(r^{-(n+1)})$ and
     $\int_{D_r}  |\gra^2 Rc|^2   dV= O(r^{-(n+4)})$ we have $|Rm| = O(r^{-n}).$
  
   By (\ref{i:poh-main}) and Lemma~\ref{l:wKato1}, 
    \begin{align}
 &\int_{D_r}  \frac{n-2}{2} |\gra Rm|^2   dV \leq       \int_{D_r} C(|\gra Rm|^2 |\gC|  |X| 
    +  |\gra Rm|   |Rm|^2  |X| + |\hess Rc| |X| |\gra Rm|) dV  \notag    \\
  &+ \left( \frac{n-1}{n+1}- \frac{1}{2}+ Cr^{-\gd_0} \right) r \int_{S_r}   |\gra Rm|^2 d\gs +   r \int_{S_r} C  (|\gra Rc|^2 + |\gra Rc||\gra Rm|) d\gs. \notag
  \end{align} 
   Therefore, 
   \begin{align}
  &\int_{D_r} \left( \frac{n-2}{2} - \epsilon \right) |\gra Rm|^2   dV 
   \leq     \int_{D_r} C(| Rm|^4   |X|^2   +   |\hess Rc|^2 |X|^2 ) dV  \notag\\
     &+   \left(\frac{n-1}{n+1}- \frac{1}{2} + \epsilon \right) r \int_{S_r}   |\gra Rm|^2 d\gs +  r \int_{S_r}  C |\gra Rc|^2 d\gs  
        \label{i:poh-2}
      \end{align} for some $\gep > 0$ small.
   Hence,  for some $\gep' > 0$ small.        
  $$  \int_{D_r} |\gra Rm|^2 dV \leq   
  \left(\frac{(n+1)(n-2)}{n-3} - \epsilon' \right)^{-1} r \int_{S_r}   |\gra Rm|^2 d\gs+
      O( r^{-(n+2)}+r^{-(4 \ga-2-n)}).$$
   Since $4 \ga-2-n= 3n -2 -\tilde \gep$ and $\frac{(n+1)(n-2)}{n-3} - \epsilon'$ are both strictly larger than  $n+2,$ 
   by Lemma~\ref{l:ode}  we have $\int_{D_r}  |\gra Rm|^2 dV= O(r^{-(n+2)}).$ Finally, by (\ref{i:Rm}), $|Rm|= O(r^{-n}).$
 \epf
 
  \bpf[Proof of Theorem~\ref{t:Kahler}]   
    We will use the coupled system (\ref{e:Rc}).
    By (\ref{CW}) and the work by Tian-Viaclovsky \cite{TV06}, we first get that $(M, g)$ has 
    the maximum volume growth and is ALE of order zero.    
    Then harmonic  metrics case follows by Theorem~\ref{t:delRm} by letting $s \rtarw \infty.$
    It remains to prove the case of Kahler metrics and (anti-)self-dual metrics.

   Since $Rc$ satisfies $\Delta Rc = Rc \ast Rm,$ we have
    $$\Delta |Rc| \geq -C |Rc| |Rm|. $$ By Lemma~\ref{l:kahler} and \ref{l:self-dual},
    $$\Delta |Rc|^{1- \frac{2}{n}} \geq -C |Rc|^{1- \frac{2}{n}} |Rm|.$$
       By Lemma~\ref{l:homo} (b),
    we have $|Rc|^{1- \frac{2}{n}} = O(r^{-\ga_0})$ for all $\ga_0 < n-2.$
   Thus $|Rc|= O(r^{-\ga_1})$ for all $\ga_1 < n.$ 
   By Lemma~\ref{l:homo} (a) and  the Sobolev inequality. 
    \beq\label{i:Rc}
 \sup_{D_{2r}} |Rc| \leq Cr^{-\frac{n-2}{2}} \|Rc\|_{L^\frac{2n}{n-2}(D_r)} \leq Cr^{-\frac{n-2}{2}} \|\gra Rc\|_{L^2(D_r)}. 
 \eeq 
   Let $\gc= \frac{n}{n-2} > 1.$
    
    Case 1. $n \geq 5.$
   
   Applying Lemma~\ref{l:inhomo} (c) to the equation
   $\Delta |Rm| \leq-C |Rm|^2 -C |\hess Rc|$ and using (\ref{Rc-reg}) with $q= \frac{n}{2} \gc ,$ 
   we get that $|Rm| = O(r^{-\ga_2})$ for all $\ga_2 < n-2.$
   Therefore, $\ga_2 > 2.$
   By (\ref{CW}), we have $|\gra^k Rm|= O(r^{-(\ga_2 +k)}).$ Hence, by \cite{BKN89}, 
    there exist coordinates $x$ and 
    $0<\gd_0<1$ such that  $C^{-1} r \leq |x| \leq C r,$ $|g-\gd| = O(r^{-\gd_0})$ and $|\de g| = O(r^{-\gd_0 -1}).$ 
   Since distances $|x|$ and $r$ are equivalent, we will still denote $|x|$ by $r$ and denote
     by $D_r$ the complement of the coordinates ball of radius $r$ and by $S_r$
     the coordinates sphere of radius $r.$

       Applying Proposition~\ref{p:poh} with
     $T= Rc$ and $\Omega= D_r$ and using (\ref{e:Rc}), we obtain        
     \beqn 
   & &\int_{D_r}  Rc \ast Rm \ast X \ast \gra Rc \, dV= \int_{D_r}  \Delta R_{ij} X_k R_{ij,k}  dV\\
   &\geq& \int_{D_r} \left( \frac{n-2}{2} |\gra Rc|^2 - C|\gra Rc|^2 |\gC|  |X| -  C|\gra Rc|   |Rc|  |Rm|  |X| \right) dV\\
  &+& \int_{S_r}   \nu_l R_{ij, l} X_k R_{ij,k} d\gs-\int_{S_r} \frac{1}{2}  X_k \nu_k  |\gra Rc|^2 d\gs,
 \eeqn where $\nu$ is the unit outward normal.  Note that  $\nu= - \frac{X}{r}+ O(r^{- \gd_0}).$ 
 By Lemma~\ref{l:kahler} and \ref{l:self-dual} again, the above formula becomes
   \beqn 
   \int_{D_r}  \frac{n-2}{2} |\gra Rc|^2   dV &\leq& 
     \left( \frac{n}{n+2}- \frac{1}{2} \right) r\int_{S_r}   |\gra Rc|^2 d\gs+ C\int_{S_r} |\gra Rc|^2 r^{1- \gd_0} d\gs\\
    &+& C \int_{D_r} \left( |\gra Rc|^2 |\gC|  |X| +  |\gra Rc|   |Rc|  |Rm|  |X|\right) dV.
 \eeqn 
  Since $|x|$ and the geodesic distance are equivalent, 
   we can apply (\ref{Rc-reg}), (\ref{Rm-reg}) and (\ref{CW}) on $D_r,$ the complement of coordinates ball. 
  Use (\ref{Rc-reg}) to get that $\|\gra Rc\|_{L^2 (D_{2r})} \leq$ \\$ C r^{-\frac{2}{\gc}+ \frac{n}{2}}  \|\gra Rc\|_{L^{\frac{n}{2} \gc}(D_{2r})}
  \leq C r^{-3 + \frac{n}{2}} \|Rc\|_{L^{\frac{n}{2}}(D_r)}= O(r^{-(\ga_1 + 1 - \frac{n}{2})}).$ 
   Observe that $\gra Rc$ satisfies $\Delta \gra Rc = Rm \ast \gra Rc + \gra Rm \ast Rc.$
 Thus,  $\Delta |\gra Rc| \geq -C |Rm||\gra Rc| -C |\gra Rm \ast Rc|.$
 Applying Lemma~\ref{l:inhomo} (a) to the above equation with $p=2$ and $q= \frac{n}{2} \gc,$
   and using $\sup_{|x|= r} |\gra Rm|= O(r^{-(\ga_2 +1)})$  we get 
   \begin{align} 
   \sup_{D_{2r}} |\gra Rc| &\leq Cr^{-{\frac{n}{2}}} \|\gra Rc\|_{L^2(D_r)} +C r^2 \sup_{D_r} |\gra Rm||Rc| \notag \\
   &= O(r^{-(\ga_1+1)}+r^{-(\ga_2+\ga_1-1)})= O(r^{-(\ga_1+1)}).\label{i:Rc-reg1}
   \end{align}
  Therefore,
 $$  \int_{D_r}   |\gra Rc|^2   dV \leq \frac{1}{n+2} r \int_{S_r}   |\gra Rc|^2 d\gs+ O(r^{-(2 \ga_1+2+ \gd_0-n}+ r^{-(2 \ga_1+ \ga_2-n)}).
 $$
 Since   $2 \ga_1+ \ga_2-n> 2 \ga_1+2+ \gd_0-n> n+2,$ by Lemma~\ref{l:ode} with 
 $f = \int_{D_r} |\gra Rc|^2 dV,$ we have 
  $\int_{D_r}  |\gra Rc|^2 dV= O(r^{-(n+2)})$ and by (\ref{i:Rc}), $|Rc|= O(r^{-n}).$ (Since $|x|$ and the geodesic distance are equivalent,
   we can apply  (\ref{i:Rc}) on $D_r,$ the complement of coordinates ball.)
 Note that  the above $Rc$ estimates are independent of the dimensions once we have 
 $|g - \gd| = O(r^{-\gd_0}).$

 Now we prove that $|Rm|= O(r^{-n}).$ Similar as before, we have $\sup_{|x| = r} |\gra Rc| = O(r^{-(n+1)})$ by (\ref{i:Rc-reg1}) with $\ga_1 = n$ now.
    By (\ref{Rc-reg}), we also have
    \beq \label{i:Rc-reg2}
    \|\hess Rc\|_{L^2(D_{2r})} \leq Cr^{-\frac{2}{\gc}+ \frac{n}{2}}  \|\hess Rc\|_{L^{\frac{n}{2} \gc}(D_{2r})}
  \leq C r^{-4 + \frac{n}{2}} \|Rc\|_{L^{\frac{n}{2}}(D_r)}= O(r^{-(\frac{n}{2}+ 2)}).
    \eeq
   By (\ref{i:poh-2}), we obtain
   \beqn
   \int_{D_r} \left( \frac{n-2}{2} - \epsilon \right) |\gra Rm|^2   dV 
   &\leq&     \left(\frac{n-1}{n+1}- \frac{1}{2} + \epsilon \right) r \int_{S_r}   |\gra Rm|^2 d\gs \\
   &+&  O(r^{-(n+2)}+ r^{-(4 \ga_2- 2- n)}).
    \eeqn
    Therefore,
    \beq \label{i:ode3}
    \int_{D_r}  |\gra Rm|^2   dV 
   \leq     \left(\frac{(n+1)(n-2)}{n-3}- \epsilon' \right)^{-1} r \int_{S_r}   |\gra Rm|^2 d\gs 
   +  O(r^{-(n+2)}+ r^{-(4 \ga_2- 2- n)})
    \eeq for $\gep'$ small. Note that $\frac{(n+1)(n-2)}{n-3}- \epsilon' > n+2.$

   When $n \geq 7,$ $4 \ga_2- 2- n= 3n -10 - \tilde \gep > n+2.$ 
    By Lemma~\ref{l:ode}, we get 
  $\int_{D_r}  |\gra Rm|^2  = O(r^{- (n+2)}).$ 
    Applying Lemma~\ref{l:inhomo} (a) to $\Delta |Rm| \geq -C |Rm|^2 -C |\hess Rc|$ with $p = \frac{2n}{n-2}$ and $q= \frac{n}{2} \gc,$ we have
 \beq\label{i:Rm'}
 \sup_{D_{2r}} |Rm| \leq Cr^{-\frac{n-2}{2}} \|Rm\|_{L^\frac{2n}{n-2}(D_r)} + C r^{-n} \leq Cr^{-\frac{n-2}{2}} \|\gra Rm\|_{L^2(D_r)} + C r^{-n}. 
 \eeq Therefore, $|Rm|= O(r^{-n}).$
  
   When $n=5$ and $6,$  $4 \ga_2- 2- n= 3n -10 - \tilde \gep < n+2.$ 
   By Lemma~\ref{l:ode} and (\ref{i:Rm'}), $\int_{D_r}  |\gra Rm|^2  = O(r^{- (3n -10 - \tilde \gep)})$ and $|Rm|= O(r^{-(2n-6-\tilde \gep)}).$
   Now going back to (\ref{i:ode3}) with $\ga_2 = 2n-6-\tilde \gep,$ we have  $4 \ga_2- 2- n= 7n -26 - 4\tilde \gep > n+2.$
   Thus, by Lemma~\ref{l:ode} and (\ref{i:Rm'})  $|Rm|= O(r^{-n}).$   
  
 Case 2. $n=4.$
  
  Consider the equation  $\Delta |Rm| \geq -C |Rm|^2 -C |\hess Rc|.$
   By (\ref{i:Rc-reg2}),
   $ \|\hess Rc\|_{L^2(D_{2r})} \leq\\ C r^{-2} \|Rc\|_{L^{\frac{n}{2}}(D_r)}= O(r^{-\ga_1}).
    $ Therefore,
   by Lemma~\ref{l:inhomo} (a) and (\ref{Rm-reg}), we have
  \begin{align}
  \sup_{D_{2r}} |Rm| \leq Cr^{-4/p} \|Rm\|_{L^p(D_r)} + Cr^{-\ga_1},\notag \\
  \|\gra Rm\|_{L^2 (B_{r/2} (x))}  \leq C r\|\gra Rm\|_{L^4(B_{r/2} (x))} \leq C r^{-1}  \|Rm\|_{L^2(B_r (x))},  \label{i:dim4'}
  \end{align} for $p> 1.$
  Since $|Rm|^2$ is integrable, let $\int_{D_r} |Rm|^2 dV= \gep_0(r)^2,$ where $\gep_0 (r) \rtarw 0$ when $r \rtarw \infty.$
  Therefore, by \cite{TV06} we have $|g- \gd| \leq C \gep_0 (r), |\de g| \leq C \gep_0 (r)/r$ and $ |Rm| \leq C \gep_0 (r)/r^2.$

  By (\ref{i:poh-dim4}),
    \begin{align} 
 &\int_{D_r}   |\gra Rm|^2   dV \leq  \int_{D_r} C(|\gra Rm|^2 |\gC|  |X| 
    +  |\gra Rm|   |Rm|^2  |X| + |\gra \gd Rm| |X| |\gra Rm|) dV  \notag    \\
  & + (\frac{1}{6}+ C\epsilon_0 (r)) r \int_{S_r}   |\gra Rm|^2 d\gs +   r \int_{S_r} C  (|\gd Rm|^2 + |\gd Rm||\gra Rm|) d\gs. \label{i:poh-dim4'}
   \end{align}
  Let $A_r = \{x: r/2 < |x| < 3r/2\}.$ By (\ref{i:dim4'}) and Sobolev inequality,
  \beqn
  \int_{A_r}  |\gra Rm|   |Rm|^2 r dV &\leq& C r \|\gra Rm\|_{L^2(A_r)} \|Rm\|_{L^4(A_r)}^2 \leq C r \|\gra Rm\|_{L^2(A_r)}^3\\
   &\leq & C \|Rm\|_{L^2(B_{\frac{5}{4} r} \setminus B_{\frac{3}{4} r})}   \|\gra Rm\|_{L^2(A_r)}^2. 
  \eeqn
  Hence,  $$ \int_{D_r}  |\gra Rm|   |Rm|^2 |X| dV \leq  C  \gep_0 (r)   \|\gra Rm\|_{L^2(D_r)}^2.$$
  Now (\ref{i:poh-dim4'})  becomes
  \begin{align} 
  (1-\gep_0 (r)- \gep)  \int_{D_r}   |\gra Rm|^2  dV & \leq  C \int_{D_r} |\gra \gd Rm|^2 |X|^2 dV\notag \\
   & + (\gep+ \frac{1}{6}) r \int_{S_r}   |\gra Rm|^2 d\gs 
  +   r \int_{S_r} C  |\gd Rm|^2 d\gs \label{i:poh-dim4''}
  \end{align} for some $\gep > 0$ small.
  To compute the above formula,
   by (\ref{Rc-reg}), 
    $\|\gra Rc\|_{L^2(D_{2r})} \leq$ \\$ Cr  \|\gra Rc\|_{L^4 (D_{2r})}
  \leq C r^{-1} \|Rc\|_{L^2(D_r)}= O(r^{-(\ga_1- 1)}).$
   By (\ref{i:Rc-reg1}) and $|\gra Rm| \leq C \epsilon_0 (r)/ r^3,$  we get 
   \begin{align}  
   \sup_{D_{2r}} |\gra Rc| &\leq Cr^{-2} \|\gra Rc\|_{L^2(D_r)} +C r^2 \sup_{D_r} |\gra Rm| |Rc| \notag \\
   &= O((1+ \gep_0(r))r^{-(\ga_1+1))})= O(r^{-(\ga_1+1)}). \label{i:Rc-reg1''}
   \end{align}
  Go back to (\ref{i:poh-dim4''}). Noting that $|\gd Rm| \leq |\gra Rc|,$ $|\gra \gd Rm|\leq |\hess Rc|$ and using (\ref{i:Rc-reg1''})
  and   $ \|\hess Rc\|_{L^2(D_r)} = O(r^{-\ga_1})$ we obtain
    $$
    \int_{D_r}  |\gra Rm|^2   dV 
 \leq (6- \epsilon')^{-1}   r \int_{S_r}   |\gra Rm|^2 d\gs   +   O( r^{-(2\ga_1 -2)}) 
    $$ for $\epsilon'$ small.  
  Since  $2\ga_1 -2= 6 - \tilde \gep > 6- \epsilon'$ by choosing $\tilde \gep \ll \epsilon',$
  by Lemma~\ref{l:ode} and  (\ref{i:Rm'}), $|Rm|= O(r^{-(4 - \frac{\epsilon'}{2})}).$ Hence, $|Rm|= O(r^{-\ga})$ for all $\ga < 4.$
  
  By \cite{BKN89}, 
    there exist coordinates $x$ and  $0<\gd_0<1$ such that $C^{-1} r \leq |x| \leq C r,$ $|g-\gd| = O(r^{-\gd_0})$ and $|\de g| = O(r^{-\gd_0 -1}).$ 
 Now by the same argument  as in Case 1 ($n \geq 5$), (the $Rc$ estimates are independent of the dimensions once we have 
 $|g - \gd| = O(r^{-\gd_0})$)  we get $|Rc|= O(r^{-4}),$ $|\gra Rc|= O(r^{-5})$ and $\|\hess Rc\|_{L^2} = O(r^{-4}).$ 
 Then by  the same proof as in Theorem~\ref{t:delRm} (b), we get $|Rm| = O(r^{-4}).$
 
 Finally, for both Cases 1 and 2 by \cite{BKN89} there exist coordinates of order $n-2.$
 \epf

\section{Proof of Corollary~\ref{c:remov}} \label{s:cor}
  \bpf[Proof of Corollary~\ref{c:remov}]
  By going to the universal cover, we may assume the group $\Gamma= \{e\}.$
 Since $g$ is of $C^0,$ in particular $\det g$ is of $C^0$ and thus $Vol B_r \leq C r^n.$ 
  Choose $r$ small such that $(\int_{B_r} |Rm|^{n/2} dV)^{2/n} < \epsilon_0,$ where $\epsilon_0$ is as
   in Lemma~\ref{l:sing-homo} and \ref{l:sing-inhomo}.
   Applying Lemma~\ref{l:sing-homo} to the equation 
   $$ \Delta |Rc| \geq - C |Rm| |Rc|$$
   gives $|Rc| = O(r^{-\ga})$ for all $\ga > 0.$ 
   Let $A_r = \{x: r/2< |x| < 3r/2\}.$ Let $q > n/2.$ By (\ref{Rc-reg}) 
   $$\|\hess Rc\|_{L^q (A_r)} \leq C r^{-4 + \frac{n}{q}} \|Rc\|_{L^{\frac{n}{2}} (B_{2r} \setminus B_{\frac{r}{4}})} \leq C r^{-2 - \ga + \frac{n}{q}}$$
   for all $\ga > 0.$
   Now applying Lemma~\ref{l:sing-inhomo} to the equation
   $$ \Delta |Rm| \geq - C |Rm| |Rm| -C |\hess Rc|$$ 
   produces $|Rm| = O(r^{- \ga})$ for all $\ga > 0.$ Therefore, by \cite{BKN89} P342 (replacing $|x|^2$ by $|x|^{2- \ga}$ in our case)
   there  exists a diffeomorphism $\phi$
   of $B_r$ such that $\phi^{\ast} g - \gd = O(r^{-\ga +2})$ and $\de\,  \phi^{\ast} g = O(r^{-\ga +1}).$
   Hence, $g \in C^{1, \gb}$ for some $0 <\gb< 1.$
    By DeTurck-Kazdan \cite{DK81}, there exist  harmonic coordinates  around the origin. Now apply standard regularity 
    to
      $$ \begin{array} {l}
        \Delta g = -2 Rc + Q(g,  \de g)\\
        \Delta Rc = Rm \ast Rc
        \end{array}$$ for $g$ in harmonic coordinates.
        We first have that $Rc \in L^p(B_1)$ for all $p.$ By the second equation, $Rc \in W^{2, p}(B_1).$ 
        Going back to the first equation, the right hand side is in $W^{1, p}.$ Thus, $g \in W^{3, p}(B_1).$ Bootstrapping in this manner, we finally get that
       $g$  extends to a smooth Riemannian metric across the origin.
  \epf

 Department of Mathematics, University of California, Berkeley, CA
 \par
 Email address: \textsf{sophie@math.berkeley.edu}
 
 Current address:
 
  Institute for Advanced Study, Princeton, NJ
  \par
  Email address:  \textsf{sophie@math.ias.edu}

\end{document}